\documentclass[a4paper,Haag duality]{mathscan}
\newtheorem{thm}{Theorem}[section] 
\newtheorem{pro}[thm]{Proposition}  
\newtheorem{cor}[thm]{Corollary}    
\newtheorem{lem}[thm]{Lemma}        
\theoremstyle{definition}           

\newcommand{\NI}{\noindent}

\newcommand{\bea}{\begin{eqnarray}}
\newcommand{\eea}{\end{eqnarray}}

\def \b #1 {\bf #1}
\newcommand{\IR}{\mathbb{R}}
\newcommand{\IN}{\mathbb{N}}
\newcommand{\IE}{\mathbb{E}}
\newcommand{\IC}{\mathbb{C}}

\newcommand{\IZ}{\mathbb{Z}}
\newcommand{\cal}{\mathcal}
\newcommand{\cla}{{\cal A}}
\newcommand{\clk}{{\cal K}}

\newcommand{\clm}{{\cal M}}
\newcommand{\clz}{{\cal Z}}
\newcommand{\cli}{{\cal I}}
\newcommand{\cls}{{\cal S}}

\newcommand{\clh}{{\cal H}}
\newcommand{\clp}{{\cal P}}
\newcommand{\clo}{{\cal O}}
\newcommand{\clb}{{\cal B}}

\newcommand{\clj}{{\cal J}}
\newcommand{\cln}{{\cal N}}

\newcommand{\clc}{{\cal C}}

\newcommand{\raro}{\rightarrow}

\newcommand{\vsp}{\vskip 1em}

\newcommand{\ul}{\underline}

\def \qed {\hfill \vrule height6pt width 6pt depth 0pt}
\newcommand{\be}{\begin{equation}}
\newcommand{\ee}{\end{equation}}
\newcommand{\ben}{\begin{eqnarray*}}
\newcommand{\een}{\end{eqnarray*}}
\begin{document}

\title{ Extremal unital completely positive maps and its symmetries } 

\author{ Anilesh Mohari }
\thanks{...}

\address{ The Institute of Mathematical Sciences, CIT Campus, Taramani, Chennai-600113 }

\email{anilesh@imsc.res.in}



\thanks{  }  

\begin{abstract}

We consider the convex set of ( unital ) positive ( completely ) maps from a $C^*$ algebra $\cla$ to a von-Neumann sub-algebra $\clm$ of $\clb(\clh)$, the algebra of bounded linear operators on a Hilbert space $\clh$ and study its extreme points via its canonical 
lifting to the convex set of ( unital ) positive ( complete ) normal maps from $\hat{\cla}$ to $\clm$, where $\hat{\cla}$ is the universal enveloping von-Neumann algebra over $\cla$. If $\cla=\clm$ and a ( complete ) positive operator $\tau$ is a unique sum of a normal and a singular ( complete ) positive maps. Furthermore, a unital complete positive map is a unique convex combination of unital normal and singular complete positive maps. We used a duality argument to find a criteria for extremal elements in the convex set of unital completely positive maps having a given faithful normal invariant state. In our investigation, gauge symmetry in Stinespring representation and Kadison theorem on order isomorphism played an important role. 

\end{abstract}

\maketitle 
\section{ Introduction: } 

\vsp 
An irreversible quantum dynamics in discrete time is governed by a complete positive map [St] on a non commutative algebra of observables. A given physics problem often fixes the algebra of observables made of elements from a suitable unital $C^*$-algebra $\cla$ over the field of complex numbers. Once a unital $C^*$-algebra $\cla$ is fixed, a mathematically challenging problem is to classify or parametrize all complete positive maps by more elementary mathematical objects. E. St\o rmer [St\o] and W. Arveson [Ar] 
gave mathematical foundation to this classification problem within the frame of Krein-Milman-Choquet theory [Ph]. A simplified presentation of St\o rmer-Arveson's approach towards the same problem is also given by M. D. Choi [Ch] for $\cla=\!M_d(\IC)$, the algebra of $d$-dimensional matrices over the field of complex number $\!C$. In other similar context with possible extensions of these results, extreme points are studied 
in details in many follow up works, here we cite some of them [Pas], [HMP], [Ra], [Ts] 
for a historical account.       

\vsp 
On the other hand, a state of a bipartite system, either in quantum information theory 
[MW,Pa,PSa,Oh] or two sided infinite quantum spin chain models [Mo2,Mo3] on a lattice, give rises unital completely positive map as follows: Let $\cla_L$ and $\cla_R$ be two $C^*$-algebras of a unital $C^*$ algebra $\cla$ such that $\cla = \cla_L \otimes \cla_R$, where $C^*$ norm of $\cla$ is taken to be the maximal one. For two faithful states $\omega_L$ and $\omega_R$ on $\cla_L$ and $\cla_R$ respectively, we consider the non empty convex set 
$$C_{\omega_L,\omega_R} = \{ \omega: \cla \raro \IC, \mbox{ a state with }\;\omega_|\cla_L=\omega_L,\;\omega_|\cla_R= \omega_R \}$$    
Then each element $\omega \in C_{\omega_L,\omega_R}$ determines unique unital completely positive maps 
$$\tau_{\omega}:\cla_R \raro \pi_{\omega_L}(\cla_L)''$$
$$\tilde{\tau}_{\omega}:\cla_L \raro \pi_{\omega_R}(\cla_R)''$$
such that 
$$\omega_L \tau_{\omega} = \omega_R$$
and 
$$\omega_R \tilde{\tau}_{\omega} = \omega_L$$
Furthermore, the maps $\omega \raro \tau_{\omega},\;\omega \raro \tilde{\tau}_{\omega}$ 
are affine continuous maps in Bounded weak topology [Ar,Pa]. The map $\tau_{\omega} \raro \tilde{\tau}_{\omega}$ is an affine one to one map satisfying the duality relation for all $x \in \cla_L$ and $y \in \cla_R$
\be 
\omega(x \otimes y)= <\clj_L \tau_{\omega}(y) \clj_L \zeta_{\omega_L}, \pi_{\omega_L}(x) \zeta_{\omega_L}> = 
<\clj_R \pi_{\omega_R}(x) \clj_R \zeta_{\omega_R}, \tilde{\tau}(y) \zeta_{\omega_R}>
\ee
where $\omega_L$ and $\omega_R$ are cyclic and separating vectors for $\pi_{\omega}(\cla_L)''$ and $\pi_{\omega_R}(\cla_R)''$ in $\clh_{\omega_L}$ and $\clh_{\omega_R}$ respectively. Furthermore, $\clj_L$ and $\clj_R$ are conjugate operators of Tomita [BR1] associated with $\omega_L$ and $\omega_R$ respectively. Thus the convex set $CP_1(\cla,\clm)$ of unital completely positive maps from a $C^*$ algebra $\cla$ to a von-Neumann algebra $\clm$ plays an important role in quantum communication channels [Par,PSa,MW,Oh]. Our main aim in this paper is to find finer structures of its extreme points and investigate how two extreme points are related for a possible classification of its extreme points upto cocycle conjugation.     
    
\vsp 
In this paper, we start with a little different approach then [St\o,Ar] and generalize classical work of J. Tomiyama [To1,To2,To3] on Lebesgue decomposition of a norm one projection on a von-Neumann algebra $\clm$. We follow closely S. Sakai's Lebesgue decomposition [Ta] of a positive functional on a von Neumann algebra to prove: an extremal unital completely positive map on a von Neumann factor $\clm$ is either a normal map or a singular map. Furthermore, any other unital completely positive map $\tau$ on a von Neumann factor $\clm$ is a unique convex combination of a normal $\tau_{\sigma}$ and a singular $\tau_s$ unital completely positive map on $\clm$. In this analysis, Stinespring representation played no role. 

\vsp 
More generally, let $CP_1(\cla,\clm)$ be the convex set of unital completely positive map from a $C^*$-algebra $\cla$ to a von Neumann algebra $\clm$, which is acting on separable Hilbert space $\clh$ over the field of complex numbers $\IC$ and $\clb(\clh)$ be the algebra of bounded operators on $\clh$.

\vsp 
For a given unital completely positive map $\tau:\cla \raro \clm$, we have unique upto isomorphism minimal Stinespring representation $(\clk,\pi,V)$
\be 
\tau(x) = V \pi(x) V^*
\ee 
for all $x \in \cla$, where $\pi:\cla \raro \clb(\clk)$ is 
a unital $*$-homomorphism and $\clk$ is the cyclic Hilbert space generated by $\pi$ and
$\clh$ i.e. $\{ \pi(x)f: x \in \cla,f \in \clh \}$ is total in $\clk$ and 
$V^*:\clh \raro \clk$ is an isometry. Let $\hat{\pi}:\cla^{**} \raro \cla$ be the universal lifting map so that $\tilde{\pi} i = \pi$ on $\cla$, where $\cla^{**}$ is the universal von-Neumann algebra over $\cla$ isomorphic to the dual Banach space of 
the dual space $\cla^*$ of $\cla$. Let $z_{\pi}$ be the element in the centre of $\cla^{**}$, the support projection of the representation $\pi$, then we define elements $\tau_{\sigma},\tau_s \in CP(\cla,\clm)$ by  
$$\tau_{\sigma}(x) = V \tilde{\pi}(i(x)z_{\pi})V^*$$     
$$\tau_s(x) = V \tilde{\pi}(i(x)z_{\pi})V^*$$
for all $x \in \cla$. 
Furthermore, we prove that the collection 
$$CP^{\sigma}_1(\cla,\clm)= \{ \tau_{\sigma}: \tau \in CP_1(\cla,\clm) \}$$ 
is also a sequentially closed convex set in Bounded Weak topology and any 
element $\tau \in CP_1(\cla,\clm)$ is a convex combination of elements 
from $CP^{\sigma}_1(\cla,\clm)$ and $CP^s_1(\cla,\clm)$, where 
$$CP^s_1(\cls,\clm) = \{ \tau \in CP_1(\cla,\clm): \tau_{\sigma} = 0 \}$$  

\vsp 
In particular, using duality relation (1) and symmetries in Stinespring representation (2), we prove $\tau$ is an extremal element in 
$$CP_{\phi} = \{ \tau \in CP_1(\cla,\clm): \phi \tau = \phi \}$$
where $\clm = \pi_{\phi}(\cla)''$ for a faithful normal state of $\clm$ if, and only if there exists no non-trivial $\lambda=(\lambda^k_j)$ with entries in $\clm'$ satisfying 
the relation 
$$
\sum_{\alpha,\beta} v_{\alpha} \lambda^{\alpha}_{\beta} v^*_{\beta}=0,\;\;\sum_{\alpha,\beta} \tilde{v}_{\alpha} \tilde{\lambda}^{\alpha}_{\beta} \tilde{v}^*_{\beta}=0,
$$
where $\lambda \raro \tilde{\lambda}$ is an unital order-isomorphism map 
on $\!M_{n_{\tau}}(\clm'))$, $(v^*_{\alpha}:1 \le \alpha \le n_{\tau})$ and $(\tilde{v}^*_{\alpha}:1 \le \alpha \le n_{\tau})$ are Krause elements in 
minimal Stinespring representation of $\tau$ and $\tilde{\tau}$ respectively. 
Furthermore, if $\tau$ admits inner representation, then same condition 
holds with entries $t^{\alpha}_{\beta}$ in the centre of $\clm$ and order 
isomorphism is given by
$$\tilde{t}^{\alpha}_{\beta}=t^{\beta}_{\alpha}$$
The proof relies on St\o rmer-Arveson version of Radon-Nikodym theorem for completely positive maps and Kadison theorem [Ka1] on unital order isomorphism on $C^*$-algebras. In particular, this result generalize a well known criteria [LS] valid for $\clm=\!M_n(\IC)$ and $\phi=tr$, the normalize trace.

\vsp 
In the last section, we deal with simplest situation $C_{\phi}=C_{\omega_L,\omega_R}$
where we have taken $\cla_L=\!M_n(\IC)$ and $\cla_R=\!M_n(\IC)$ with $\omega_L=\omega_R=\phi$, where $\phi$ is a faithful state on $\!M_n(\IC)$. 
We prove that an element $\psi \in C_{\phi}$ is a factor state if $\tau_{\psi}$ is an extremal element in $CP_{\phi}$. Furthermore, $\psi$ is a pure state in $C_{\phi}$ if, and only if $\tau_{\psi} \in CP_{\phi}$ is an automorphism on $\clm$. The last statement is a generalization of a theorem proved in [Oh] with $\phi=tr$.

\section{ A decomposition of completely positive maps: } 
  
\vsp 
Let $\cla$ be a unital $C^*$-algebra. A linear map $\tau: \cla \raro \cla$ is called {\it positive } if $\tau(x) \ge 0$ for all $x \ge 0$. Such a map is automatically bounded with norm $||\tau||=||\tau(I)||.$ A linear map $\tau:\cla \raro \cla$ is called $n$-{\it positive } [St] $(CP)$ if $\tau \otimes I_n: \cla \otimes M_n \raro \cla \otimes M_n$ is positive, where $\tau \otimes I_n$ is defined by $(x^i_j) \raro ( \tau(x^i_j) )$ with elements $(x^i_j)$ are in $\cla$. If $\tau$ is $n$-positive for each $n \ge 1$, $\tau$ is called {\it completely positive}. 
 
\vsp 
For each unital representation $\pi:\cla \raro \clb(\clh_{\pi})$ of $\cla$, we denote by $\clm(\pi)$ the von-Neumann algebra $\pi(\cla)''$ generated by $\pi(\cla)$ in $\clb(\clh_{\pi})$. A representation $(\pi,\clh_{\pi})$ is called {\it universal} if for any representation $(\rho,\clk_{\rho})$ there exists a $\sigma$-weak continuous $*$-homomorphism $\hat{\rho}$ of $\clm(\pi)$ onto $\clm(\rho)$ such that 
$$\rho(x) = \hat{\rho}(\pi(x))$$ 
If $(\pi,\clh_{\pi})$ is a universal representation of $\cla$, $\clm(\pi)$ is called universal enveloping von-Neumann algebra of $\cla$. The universal enveloping von-Neumann algebra is uniquely determined upto isomorphism. We recall now two standard constructions of universal enveloping von-Neumann algebras, each one has its own merit. 

\vsp 
Let $\cla^*$ be the Banach space dual of $\cla$ and $\cla^{**}$ be the double dual of $\cla$ i.e. dual Banach space of $\cla^*$. Let $i:\cla \raro \cla^{**} \equiv \clm^{\cla}_u$ be the inclusion map of $\cla$ into  $\cla^{**}$. For a representation $\pi:\cla \raro \clb(\clh_{\pi})$, we set von-Neumann algebra $\clm_{\pi} = \pi(\cla)''$ and the Banach space adjoint linear map by
$\pi^t:\clm_{\pi}^* \raro \cla^*$. We set now linear map $(\pi^t)_*:(\clm_{\pi})_* \raro \cla^*$ by restricting $\pi^t$ 
to $(\clm_{\pi})_*.$ Finally, we set notation $\hat{\pi}:\cla^{**} \raro \clm_{\pi}$ for the Banach space adjoint map of $(\pi^t)_*.$ Thus by our construction $\hat{\pi}$ is a normal map, being the dual of a bounded linear map on their pre-dual Banach spaces.  

\vsp 
The enveloping von-Neumann algebra $\cla^{**}$ of a $C^*$-algebra $\cla$ is defined to be the double dual $\cla^{**}$ of $\cla$. That $\cla^{**}$, being a dual of a Banach space, is a von-Neumann algebra [Sa]. The following proposition says $\cla^{**}$ is indeed the
universal von-Neumann algebra of $\cla$. 

\vsp 
\begin{pro} 
Let $(\pi,\clh_{\pi})$ be a representation of $\cla$. We have the following properties:

\NI (a) $\hat{\pi}:\cla^{**} \raro \clm_{\pi}$ is a linear map which is continuous 
with respect to weak$^*$ topologies on $\cla^{**}$ and $\clm_{\pi}$. The map $\hat{\pi}$ takes the norm closed unit ball of $\cla^{**}$ onto the norm closed unit ball of $\clm_{\pi}$;

\NI (b) $\hat{\pi} \circ i = \pi$ on $\cla$;  

\NI (c) $\hat{\pi}$ is a unital completely positive map 
from $\cla^{**}$ onto $\clm_{\pi}$;

\NI (d) For any central element $z \in \cla^{**}$, $\hat{\pi}(z)$ 
is an element in the center of $\clm_{\pi}$. 
\end{pro}

\vsp 
\begin{proof} 
For statement (a) and (b), we refer to Lemma 2.2 in Chapter 3 of [Ta].

\vsp 
Statement (c) and (d) are as well known but could not find a suitable ready reference. Here we indicate a proof. Since $\pi$ is a positive map and so is its transpose $\pi^t$. Thus the restriction $(\pi^t)_*$ is also positive. That $\hat{\pi}$ is positive follows as $(\pi^t)_*$ is positive and 
onto. For $n-$positive property of $\hat{\pi}$, 
we note that the universal enveloping algebra over $M_n(\cla)$ is 
$M_n(\cla^{**})$ and the canonical map $i \otimes I_n$ is the 
inclusion map of $M_n(\cla)$ into $M_n(\cla^{**})$. Furthermore, for a 
representation $\pi:\cla \raro \clb(\clh_{\pi})$, we also have 
$\hat{\pi} \otimes I_n \circ i \otimes I_n = \pi \otimes I_n$.    
This shows in particular that $\hat{\pi}$ is a completely positive map
and its restriction on $i(\cla)$ is a representation.  

\vsp 
By Kadison-Schwarz inequality [Ka2] for unital completely positive map we have
$\hat{\pi}(x^*y)=\hat{\pi}(x^*) \hat{\pi}(y)$ for $x \in \cla^{**}$ and
$y \in i(\cla)$ since $\hat{\pi}(y^*y)=\hat{\pi}(y^*)\hat{\pi}(y)$ for $y \in 
i(\cla)$. Taking adjoint in the above relation we also get 
$\hat{\pi}(y^*x)=\hat{\pi}(y^*)\hat{\pi}(x)$. Operator algebras involved are being 
$*$-closed and also $\hat{\pi}$ being onto, we get $\hat{\pi}(x) \in \clm_{\pi} \bigcap \clm'_{\pi}$ 
if $x \in \cla^{**} \bigcap \cla^{**'}$.  
\end{proof} 

\vsp 
Let $\cls(\cla)$ be the convex set of states on $\cla$ and $(\clh_{\phi},\pi_{\phi},\zeta_{\phi})$ be the GNS representation associated with a state $\phi$ on $\cla$. The universal von-Neumann algebra of $\cla$ is given by 
$\clm^{\cla}_u = \{ \pi_u(x): x \in \cla \}'',$  
where $\pi_u(x) = \oplus_{\phi \in \cls(\cla) } \pi_{\phi}(x)$
is the direct sum of representations on
$$\clh_u = \oplus_{\phi \in \cls(\cla)} \clh_{\phi}$$
Since any representation $(\pi,\clh_{\pi})$ is a direct sum of cyclic representation, $\clm^{\cla}_u$ is the universal von-Neumann algebra of $\cla$. 

\vsp 
We defined left and right action of a $C^*$-algebra $\cla$ on its dual $\cla^*$ by 
\be 
<y, x \omega > =< yx,\omega>,\;\;\;\; <y, \omega x>=<xy,\omega>
\ee 
where $<.,.>$ denotes evaluation of a functional on a given element in the Banach space.
The crucial property that $\cla$ is an algebra is used here to define these actions on $\cla$. A subspace of $\cla^*$ is called $\cla$ {\it invariant} if the subspace is invariant by both left and right action. Given a $\cla$-invariant subspace $V$ of $\cla^*$

\vsp 
A projection in the center of $\cla^{**}$ determines uniquely a representation of $\cla$ upto quasi-equivalence i.e. a representation $\pi$ is quasi-equivalent to the sub-representation $x \raro \pi_u(x)z_{\pi}$ for some central projection $z_{\pi}$ in $\pi_{u}(\cla)''$. The projection $z_{\pi}$ is called {\it support} of the representation in $\cla^{**}$, determined unique by the representation $\pi$ upto isomorphism. Given a central projection $z_{\pi}$, there exists a uniquely determined $\cla$-invariant norm closed subspace $V(\pi)$ of $\cla^*$ given by $\pi^t(\clm_*(\pi))$ i.e. The range of the map $\pi^t: \clm^* \raro \cla^*$ restricted to $\clm_*$. Conversely, given a norm closed $\cla$-invariant subspace of $\cla^*$, there exists a central projection $z$ in $\cla^{**}$ such that the representation $\rho:x \raro \pi_{u}(x)z$ is non-degenerate and $V(\rho)=V$. For details, we once more refer to section 2 of Chapter 3 in [Ta]. 

\vsp
Let $\clm$ be a von-Neumann algebra acting on a Hilbert space $(\clh,<.,.>)$ over the field of complex numbers, where the inner product is linear in second variable. Let $\clm_*$ be the pre-dual Banach space of $\clm$. For a von-Neumann algebra $\clm$, $\clm_*$ is a proper subset unless $\clm$ is finite dimensional. However the unit ball of $\clm_*$ is a dense subset in the unit ball of $\clm^*$ in the weak$^*$ topology on $\clm^*$. Nevertheless $\clm_*$ is sequentially closed [Ta] and given any bounded linear functional $\omega$ on $\clm$ there is a unique element $\omega_{\sigma} \in \clm_*$ so that $\omega = \omega_{\sigma} + \omega_s$ determined by  
$||\omega - \omega_{\sigma}||=\mbox{min}_{\omega' \in \clm_*}||\omega-\omega'||$. 
The element $\omega_s$ is called {\it purely singular } on $\clm$ unless it is the zero element. In the following we describe such a decomposition via universal enveloping algebra of $\clm$. 

\vsp 
Now we take $\cla=\clm$, a von-Neumann algebra in Proposition 2.1 and $\pi:\clm \raro \clm 
\subseteq \clb(\clh)$ be the identity representation. Let $z_{\pi}$ be the support projection of the representation $\pi:\clm \raro \clm$ in $\clm^{**}$. Thus $\hat{\pi}(z)$ is an element in the center of $\clm$, where $\hat{\pi}$ is the lift of $\pi$ to the universal enveloping algebra defined in Proposition 2.1. 

\vsp 
\begin{pro} 
Let $\clm$ be a von-Neumann algebra with its pre-dual space $\clm_*$ and dual $\clm^*$. Then there exist 
a unique central projection $z \in \clm^{**}$ so that 
$$\clm_*= \clm^*z,$$ 
where $\clm^{**}$ acts natural on the Banach space $\clm^*$ given 
by 
$$<y, x \omega > = < yx,\omega>,\;\;\;\;\; <y, \omega x>=<xy,\omega>$$ 
for $x,y \in \clm^{**}$ and $\omega \in \clm^*$.
\end{pro} 

\vsp 
\begin{proof} 
The von-Neumann algebra $\clm$ is dense in weak$^*$ topology of $\clm^{**}$ ( since unit ball of $\clm_*$ is dense in the unit ball of $\clm^*$ ) and
$\clm_*$ is a $\clm^{**}$-invariant norm closed subspace of $\clm^*$. Thus by Theorem 2.7 in Chapter -III [Ta], we conclude that $\clm_* = \clm^* z$.  
\end{proof}  

\vsp 
\begin{pro} 
An element $\omega \in \clm^*$ determines 
uniquely an element $\omega_{\sigma} \in \clm_*$ defined by 
$$\omega_{\sigma}(x)= <\omega,z i(x)>$$ 
for all $x \in \clm$ and the map $\omega \raro \omega_{\sigma} \in \clm^*$ is positive linear contractive map on the Banach space $\clm_*$. The element $\omega_s(x)= <\omega,(1-z)i(x)>$ in $\clm^*$ is singular provided $z \neq 1$. The decomposition 
$\omega=\omega_{\sigma} + \omega_s$ is also uniquely determined and 
$||\omega||=||\omega_{\sigma}||+||\omega_s||$.  
\end{pro} 
\vsp 
\begin{proof} 
We refer chapter 3 in [Ta,Chapter 3.2] for details.
\end{proof}  

\vsp 
Let $\cla$ be a $C^*$ algebra and $\clm$ be a von-Neumann algebra acting on a Hilbert space $\clh$. For a positive map $\tau:\cla \raro \clm \subseteq \clb(\clh)$, we set a positive map $\hat{\tau}:\cla^{**} \raro \clm$, the adjoint map of 
$(\tau^t)_*:\clm_* \raro \cla^*$, where $\tau^t: \clm^* \raro \cla^*$ is the adjoint map of $\tau:\cla \raro \clm$. 
Thus we have the lifting property $$\hat{\tau} \circ i = \tau$$
where we recall $i: \cla \raro \cla^{**}$ is the canonical inclusion map.   

\vsp 
A positive linear map $\tau:\clm \raro \clm$ is called {\it normal} if $l.u.b.\tau(x_{\alpha}) = \tau(l.u.b.x_{\alpha})$ for any bounded increasing net $x_{\alpha}$ in $\clm$ where $l.u.b.$ denotes least upper bound. We will use notations $P^{\sigma}$ and 
$CP^{\sigma}(\clm)$ for the convex set of positive and completely positive normal maps 
on $\clm$. When there is no room for confusion,, we will simply denote them by $P^{\sigma}$ and $CP^{\sigma}$ respectively.  

\vsp 
In contrast, a positive non-zero map $\tau:\clm \raro \clm$ is called {\it singular} if $\phi \tau$ is either zero or a purely singular positive functional for any positive normal functional $\phi$ of $\clm$. We will use notations $P^s(\clm)$ and 
$CP^s(\clm)$ for the convex set of positive and completely positive singular maps 
on $\clm$ respectively. When there is no room for confusion, we will simply denote them by $P^s$ and $CP^s$ respectively.

\vsp 
The following proposition says along the same vein now for the class of positive or completely positive maps on $\clm$. The result is well known for quite some time [To2] for completely positive map.  

\vsp 
\begin{pro}  
Given a positive map $\tau$ on $\clm$, there is a unique normal positive map $\tau_{\sigma}$ on $\clm$ and a singular positive map $\tau_s$ such that 
$$\tau = \tau_{\sigma} + \tau_{s},$$  
where $\tau_{\sigma},\tau_{s}$ are determined uniquely by the decomposition 
$$\omega \tau_{\sigma} = (\omega \tau)_{\sigma},\;\omega \tau_s= (\omega \tau)_s$$
of a normal state $\omega$. Furthermore, the set of normal positive maps on $\clm$ is sequentially closed in the set of positive maps on $\clm$. Furthermore, if $\tau$ is 
$n$-positive, then $\tau_{\sigma}$ and $\tau_{s}$ are also $n$-positive.  
\end{pro}

\vsp 
\begin{proof} 
Let $\clm^{**}$ be the universal von-Neumann algebra of $\clm$ and 
$i:\clm \raro \clm^{**}$ be the canonical inclusion map of $\clm$ into 
$\clm^{**}$. We write for an element $\omega \in \clm_*$ 
$$(\omega \tau)_{\sigma}(x) = <\omega \tau, z i(x)>$$ 
where $z$ is the central projection in $\clm^{**}$. The map 
$\omega \raro (\omega \tau)_{\sigma}$ determines a normal 
positive linear map $\tau^{\sigma}$ on $\clm$ by 
$\omega \tau_{\sigma}(x)= <\omega \tau, z i(x)>$ 
for all $x \in \clm$ and $\omega \in \clm_*$.   
It also shows clearly that $\tau_{\sigma}$ is $n$-positive as 
$\tau_{\sigma} \otimes I_n = (\tau \otimes I_n)_{\sigma}$ for $n$-positive map $\tau$. Similarly we also have $\omega \tau_s(x) = <\omega \tau, (1-z)i(x)>$ 
for all $x \in \clm,\;\omega \in \clm_*$ and the induced map $\tau_s$ 
is also $n$-positive if $\tau$ is so.  

\vsp 
Let $\tau_n$ be a sequence of such normal maps on $\clm$ and its bounded weak limit is $\tau$ i.e. $\omega(\tau_n(x)) \raro \omega(\tau(x))$ for
all $\omega \in \clm_*$ and $x \in \clm$. Let $\tau_s$ be the singular part of $\tau$. Then singular part of $\omega \tau_n$ also converges to singular 
part of $\omega \tau$ in weak$^*$ topology. Since each $\tau_n$ is normal and so is $\omega$, we get $\omega \tau^n_s=0$. Since this holds for all 
$\omega \in \clm_*$, we arrive at our desired claim that $\tau_s=0$. 
\end{proof} 

\vsp 
A non-zero positive map is called {\it purely singular} in short {\it singular} if its normal part of the unique decomposition contributes the zero map. The convex set of unital completely positive maps on $\clm$ is compact in Bounded-Weak topology of William Arveson [Ar, Pa]. One natural question that arises at this point: when a unital completely positive map on $\clm$ is a convex combination of unital normal and singular completely positive maps? 

\vsp 
\begin{pro} 
Let $\clm$ be von-Neumann algebra and $z$ be central projection in the universal enveloping algebra 
$\clm^{**}$ such that $\clm_*=z \clm^*$. A positive map $\tau:\clm \raro \clm$ is a convex combination of two unital positive maps 
$\tau_{\sigma} \in CP_{\sigma}$ and $\tau_{s} \in CP_{s}$ with $\lambda \in [0,1]$ i.e. \be 
\tau= \lambda \tau_{\sigma} + (1-\lambda)\tau_s
\ee 
if and only if $\hat{\tau}(z)$ is a scaler where $\hat{\tau}$ is the unique completely positive unital normal map from $\clm^{**}$ to 
$\clm$ such that following relation holds:
\be
<\omega \tau, x>=<\omega, \hat{\tau}(x)>
\ee 
for all $x \in \clm^{**}$ and $\omega \in \clm_*$. Furthermore, the map $\tau \raro \hat{\tau}$ is affine and one to one. 
\end{pro} 

\vsp 
\begin{proof} 
By Proposition 2.4, we have  
$$(\omega \tau)_{\sigma}(x)= <\omega \tau, z i(x)>=<\omega, \hat{\tau}(z i(x))>$$ 
for all $x \in \clm$ and $\omega \in \clm_*$. Uniqueness of the decomposition into normal and singular completely positive maps ensures that 
$$\tau_{\sigma}(x)= \hat{\tau}(z i(x))$$ 
for $x \in \clm$. So we also have 
$$\tau_s(x) = \hat{\tau}((I-z)i(x))$$
$x \in \clm$ such that $\tau=\tau_{\sigma} + \tau_s$. 

\vsp 
If $\hat{\tau}(z) \in \clm$ 
is a scaler say $\lambda = \hat{\tau}(z) \in [0,1]$, we can rewrite $\tau = \tau_{\sigma}+\tau_s$ as in (4) by obvious modifications. On the converse, if 
$\tau$ is a convex combination of two unital elements given as in (4), we conclude 
by uniqueness of the decomposition, $\hat{\tau}(z)= \lambda$ is a scaler.  
\end{proof} 

\vsp 
A unital completely positive map $\tau:\cla \raro \clb(\clh)$ on a unital $C^*$ algebra $\cla$ admits a minimal Stinespring representation [St]  
\be 
\tau(x)= V \pi(x) V^*,
\ee
where $\pi:\cla \raro \clb(\clh^{\tau})$ is a $*$-representation of $\cla$ into $\clb(\clh^{\tau})$, which is the Hilbert space completion of the kernel given on the set $\cla \otimes \clh$ by 
$$k(x \otimes \zeta, y \otimes \eta )= <\zeta,\tau(x^*y)\eta>$$ 
and $V^*:\clh \raro \clh^s$ is an isometry from $\clh$ into $\clh^s$ defined by 
$$V^*: \zeta \raro  I \otimes \zeta$$ 
such that $\{ \pi(x) V^* \zeta: x \in \cla, \zeta \in \clh \}$ spans $\clh^{\tau}$. Such a minimal Stinespring triplet $(\clh^{\tau},\pi,V^*)$ is uniquely determined modulo unitary equivalence i.e. if $(\clh',\pi',V'^*)$ be another triplet associated with $\tau$ with cyclic property, then we get $W^*: \pi(x) V^* \zeta \raro \pi'(x)V'^*\zeta,\;x \in \cla,\;\zeta \in \clh$ is an unitary operator so that 
\be 
W^*V^*=V'^*\;\mbox{and}\;\; W^*\pi(x)W = \pi'(x),\;x \in \cla
\ee 

\vsp 
An alternative constructions of $\tau_{\sigma}$ and $\tau_s$ for unital CP map $\tau:\clm \raro \clm$ are given as follows. 

\vsp 
Let $\tau(x)=V \pi(x) V^*$ be the Stinespring minimal representation of $\tau:\clm \raro \clb(\clh)$ i.e. $\pi:\clm \raro \clb(\clk)$ is the unique 
unital representation in a Hilbert space $\clk$ and $V^*:\clk \raro \clh$ is the unique isometry ( modulo unitary equivalence ). Let $z_{\pi}$ be the support projection of the representation $\pi:\clm \raro \clb(\clk)$ i.e. is quasi-equivalent to $\pi':\clm \raro \clm^{**}$ defined by 
$$\pi'(x)= \pi_u(x)z_{\pi}$$
Note that support projection does not depend on the choice of the minimal Stinespring dilation $\pi$ we choose as $z_{\pi}=z_{\pi'}$ by Proposition 2.12 in Chapter 3 [Ta] for 
any choice $\pi'$. 

\vsp 
Furthermore, by the duality relation (5), we verify the universal lifting property of $\tau \raro \hat{\tau}$ as follows: for $\omega \in \clm_{*}$ and $x \in \clm$
$$<\omega, \hat{\tau}(i(x))>$$
$$=<\omega \tau, i(x) >$$
$$=<\omega \tau, x>$$
i.e. $\hat{\tau} \circ i = \tau$ 

\vsp 
In particular, we have by Proposition 2.5, 
$$\hat{\pi}_{\sigma}(x)=\hat{\pi}(i(x)z_{\pi})$$
$$\hat{\pi}_s(x)=\hat{\pi}(i(x)(I-z_{\pi}))$$
for all $x \in \clm$ and thus 
$$\tau_{\sigma}(x)=V \hat{\pi}(z_{\pi} i(x))V^*$$
and
$$\tau_{s}(x)=V\hat{\pi}((I-z_{\pi})i(x))V^*$$
where $z_{\pi}$ is the support of $\pi$ in $\clm^{**}$ and 
$\hat{\pi}:\clm^{**} \raro \clb(\clk)$ is the lifting map of the representation 
$\pi:\clm \raro \clb(\clk)$ defined as in Proposition 2.1. By Proposition 2.1 (d), projection $\hat{\pi}(z_{\pi})$ is an element in the center of $\clm_{\pi}=\{\pi(x):x \in \clm \}''$. So in particular, we have $\hat{\tau}(z_{\pi})=V \hat{\pi}(z_{\pi}) V^* \in \clm$. We recall by Proposition 2.1 (d), $\hat{\pi}(z_{\pi})$ is an element in the centre of $\clm$. Thus $\tau_{\sigma}(I)$ 
is scaler for a factor $\clm$ and we arrive at our main result of this section.

\vsp 
\begin{thm}  
Let $\tau$ be an unital completely positive map on a von-Neumann factor $\clm$. Then there
exists a decomposition of $\tau$ given by 
$$\tau = \lambda \tau_{\sigma} + (1-\lambda)\tau_s$$
for some $\lambda \in [0,1]$ and unital normal $\tau_{\sigma}$ and singular $\tau_s$ completely positive maps respectively. Furthermore, if $\lambda \in (0,1)$ then such a decomposition is unique.
\end{thm}

\vsp 
\section{ Radon-Nikodym theorem for completely positive unital maps: }  

\vsp 
In this section, we review a characterization of unital complete positive maps $\tau:\cla \raro \clb(\clh)$ [St\o, Ar] to be extremal elements in $CP_1(\cla,\clm)$ with some finer additional results for our main results.   

\vsp 
\begin{pro} 
Let $\tau:\cla \raro \clb(\clh)$ be a completely positive map on a $C^*$ algebra 
$\cla$ and $\eta$ be another completely positive map on $\cla$ so that for some positive constant $c$, $\eta(x) \le c \tau(x)$ for all $x \in \cla_+$. Then there exists a 
non-negative element $T' \in \pi(\cla)'$ such that 
\be 
\eta(x)=V \pi(x)T'V^*,
\ee 
where $(\clk,\pi,V^*)$ is the minimal Stinespring triplet associated with $\tau$. Furthermore, for each $x \in \clm$, $\pi(x)$ and $T'$ commutes with the representation $\rho:\tau(\cla)' \raro \clb(\clk)$ defined by
\be 
\rho(y):x \otimes f = x \otimes yf
\ee     
\end{pro} 

\vsp 
\begin{proof} 
Without loss of generality we assume that $\cla$ is a closed $*$-subalgebra of
$\clb(\clh)$, where $\clh$ is a Hilbert space. We recall Stinespring minimal representation $\tau(x)=V\pi(x)V^*$ where
$\pi:\cla \raro \clb(\clk)$ is the representation induced by the map $y \otimes f 
\raro xy \otimes f$ where $\clh$ is von-Neumann's completion of algebraic tensor product
$\cla \otimes \clh$ with the kernel
$$k(x \otimes f, y \otimes g)=<f,\tau(x^*y)g>$$
Now we set another sesqui-linear form on $\cla \otimes \clh$ defined by 
$$s(x \otimes f,y \otimes g)=<f,\eta(x^*y)g>$$
Being a kernel, we have Cauchy-Schwarz inequality to check the following steps:
$$|s(x \otimes f,y \otimes g)|^2 \le s(x \otimes f,x \otimes f)s(y \otimes g,y \otimes g) \le
c^2 ||x \otimes f|| ||y \otimes g||$$ So there exists a bounder operator $T'$ on $\clk$
such that $s(x \otimes f, y \otimes g)=<x \otimes f, T' y \otimes g>$. We claim that 
$T' \in \pi(\cla)'$. Proof follows once we check the following steps:
$$<f,\eta(x^*zy)g>=s(x \otimes f, zy \otimes g>=<x \otimes f, T' \pi(z) y \otimes g>$$
and 
$$<f,\eta(x^*zy)g>=s(z^*x \otimes f,y \otimes g>=<\pi(z^*)x \otimes f, T'y \otimes g>$$
This completes the proof.
\end{proof} 

\vsp 
We briefly recall Bounded-Weak (BW) topology [Chapter 7 in Pa] on the convex set $P_1(\cla,\clm)$ ( $CP_1(\cla,\clm)$ ) of unital positive maps ( completely map ) from a $C^*$-algebra $\cla$ to a von-Neumann algebra $\clm \subseteq \clb(\clh)$. A net of unital positive maps $\tau_{\alpha}:\cla \raro \clm \subseteq \clb(\clh)$ converges in BW topology if and only if $\tau_{\alpha}(x)$ converges to $\tau(x)$ in weak$^*$ topology of von-Neumann algebra $\clm$. Furthermore, these topology makes the convex set $P_1(\cla,\clm)$ ($CP_1(\cla,\clm)$) compact and furthermore, BW topology is metrizable if $\cla$ and as well as $\clm_*$ are separable as Banach spaces.  

\vsp 
Given an element $\tau \in P(\cla,\clm)$, there exists a unique lifting $\hat{\tau} \in P(\cla^{**},\clm)$ such that $\hat{\tau} \circ i = \tau$ on $\cla$. By our construction $\hat{\tau}$ is a normal map. Conversely, given a normal positive map 
$\hat{\tau}:\cla^{**} \raro \clm$, $\hat{\tau}$ is the unique lift of 
$\tau:\cla \raro \clm$ given by $\tau= \hat{\tau} \circ i$. 

\vsp 
Let $\tau:\cla \raro \clm$ be a completely positive map and 
$$\tau(x) = V \pi(x) V^*$$
be its Stinespring representation with $\pi:\cla \raro \clb(\clk)$ be 
a $*$-representation and $V^*:\clh \raro \clk$ be a bounded operator. 
Let $\hat{\pi}:\cla^{**} \raro \clb(\clk)$ be the universal lifting map in 
Proposition 2.1 and $z_{\pi}$ be the the support projection of the representation 
$\pi$ in $\cla^{**}$. We define two completely positive maps from 
$\cla \raro \clb(\clk)$ by    
$$\pi_{\sigma}(x) = \hat{\pi}(i(x)z_{\pi})$$   
$$\pi_s(x) = \hat{\pi}(i(x)(I-z_{\pi}))$$
for all $x \in \cla$. We set also completely positive maps from 
$\cla$ to $\clb(\clh)$ by
$$\tau_{\sigma}(x) = V \pi_{\sigma}(x) V^*$$ 
$$\tau_s(x) =V \pi_s(x) V^*$$
for all $x \in \cla$. Thus by our construction we have
\be 
\tau = \tau_{\sigma} + \tau_s
\ee

\vsp 
We say $\tau$ is {\it normal} if $\tau_s=0$. Similarly, $\tau$ is called {\it singular} if $\tau_{\sigma}=0$. The uniqueness part of the lifting theorem in Proposition 2.1 says that 
$\tau:\cla \raro \clm$ has a unique decomposition 
$$\tau = \tau_{\sigma} + \tau_s$$
as in Proposition 2.4, where $\tau_{\sigma}$ and $\tau_s$ are normal and singular CP map 
from $\cla$ to $\clb(\clh)$ defined above. Thus this notions indeed generalize the notion 
of normal and singular complete positive map to $C^*$ algebras. We use notions 
$CP^{\sigma}(\cla,\clm)$ and $CP^s(\cla,\clm)$ respectively for normal and singular CP maps. A routine now also says that $CP^{\sigma}(\cla,\clm)$ is sequentially closed in $CP(\cla,\clm)$. For a proof, we take a sequence $\tau^n \in CP^{\sigma}(\cla,\clm)$ and $\tau \in CP(\cla,\clm)$ such that $\phi \tau^n(x) \raro \phi \phi \tau(x)$ for all normal state $\phi$ on $\clm$. Since each $\phi \tau^n$ is a normal state on $\clm$, we get 
$(\phi \tau )_s=\phi \tau_s=0$. This holds for all normal states, we get $\tau_s=0$. 
     
\vsp 
\begin{pro} 
$CP^{\sigma}_1(\cla,\clm)$ is a convex face of $CP_1(\cla,\clm)$. 
\end{pro} 

\vsp 
\begin{proof} Let $\tau_0,\tau_1 \in CP^{\sigma}_1(\cla,\clm)$ and $\tau = \lambda \tau_1 + (1-\lambda)\tau_0$. We recall elements $\pi:\cla \raro \clb(\clk)$ and $V^*:\clh \raro \clk$ in Proposition 3.1 and find $\tau_0(x) = V(T_0\pi(x))V^*$
and $\tau_1(x)=V(T_1\pi(x))V^*$ for unique non-negative elements $T_0,T_1 \in \pi(\cla)'$ such that $\lambda T_1 + (1-\lambda)T_0=I$. Thus $(\tau_k)_s(x) = V(T_k \hat{\pi}(i(x)z_{\pi})V^*=0$ for $k=0,1$. By adding them up, we get $\tau_s=0$. This shows $CP^{\sigma}_1(\cla,\clm)$ is a convex set.  

\vsp 
By Proposition 3.1, if $\tau \in CP^{\sigma}(\cla,\clm)$ and $\eta \in CP(\cla,\clm)$ 
such that $\eta \le \tau$ then $\eta \in CP^{\sigma}(\cla,\clm)$. Let $\tau \in CP^{\sigma}_1(\cla,\clm)$ and $\tau_0,\tau_1 \in CP_1(\cla,\clm)$ so that
$\tau=\lambda \tau_1 +(1-\lambda)\tau_0$ for some $\lambda \in (0,1)$. Then 
we have $\tau_1 \le {1 \over \lambda}\tau$ and also $\tau_0 \le {1 \over 1-\lambda}\tau$ on non-negative elements $\cla^{**}_+$. Thus both $\tau_1,\tau_1 \in CP^{\sigma}(\cla,\clm)$. This completes the proof. 
\end{proof} 

\vsp 
As an application of Theorem 2.6 and Choquet theorem [Ph], we have the following theorem.  

\vsp 
\begin{thm} 
Let $\tau \in CP_1(\cla,\clm)$ be a unital CP map from a unital $C^*$ algebra $\cla$ to a von-Neumann algebra $\clm$ acting on a Hilbert space $\clh$. Then the following hold: 

\vsp 
\NI (a) If $\clm$ is a factor, then we have 
$$\tau = \lambda \tau_{\sigma} + (1-\lambda)\tau_s$$
for some unital normal $\tau_{\sigma}$ and singular $\tau_s$ CP maps from $\cla$ into $\clm$. Furthermore, if $\lambda \in (0,1)$, then $\tau_{\sigma}$ and $\tau_s$ are determined uniquely. 

\vsp 
\NI (b) If $\cla$ and $\clm_*$ are separable as Banach space, then there exists a regular Borel probability measure $\mu$ on the set of extreme points $CP^e_1(\cla,\clm)$ satisfying 
$$\tau = \int_{\tau_e \in CP^e_1(\cla,\clm)} \tau_e d\mu(e)$$
Furthermore, if $\tau$ is normal (respectively singular), the support of the regular measure $\mu$ is also on the set of normal (respectively singular) extremal elements. 
In particular, we also have the decomposition given in (a). 
\end{thm} 

\vsp 
\begin{proof} (a) is simple application of Theorem 2.6, where we use the uniqueness of the lifting $\tau \raro \hat{\tau}$ and affine property of the map.  First part of (b) is a simple application of Choquet theorem [Ph]. 

\vsp 
The set of unital normal completely positive maps $CP^{\sigma}_1(\cla,\clm)$ need not be a closed subset in BW topology and thus need not be compact in general in BW topology. However, for the last part, we only need to how that the set of normal extremal elements is Borel measurable. In fact, it is enough if we show that the set of unital normal maps are Borel measurable i.e. we need to show the set $\{\phi(\tau(x)): \tau \in CP^{\sigma}_1(\cla,\clm) \}$ for each $x \in \cla$ and normal state $\phi$ of $\clm$ is a Borel measurable set in $\IC$. Since the BW topology is now metrizable and the map $\tau \raro \phi(\tau(x))$ being sequentially continuous on $CP^{\sigma}_1(\cla,\clm)$, we get the map $\tau \raro \phi(\tau(x))$ is in fact continuous in BW topology and thus the set of unital normal CP maps is Borel measurable. 

\vsp 
The last part of the statement follows once we split the integral into sum of two integrals of disjoint measurable sets consists of normal and singular extreme points of $CP_1(\cla,\clm)$. \end{proof} 

\vsp 
\begin{pro} 
Let $\cla$ be a $C^*$-algebra and $\tau:\cla \raro \clb(\clh)$ be a unital completely positive map and $\tau(x)=V \pi(x) V^*$ be 
the unique upto isomorphism minimal Stinespring representation. Then following are equivalent
 
\NI (a) $\tau$ is extremal in $CP_1$; 

\NI (b) $V \Lambda V^*=0$ for $\Lambda \in \pi(\cla)'$ if and only if $\Lambda=0$;  
\end{pro}  

\vsp 
\begin{proof} 
Let $\tau=\lambda \tau_1+(1-\lambda)\tau_0$ for some $\tau_0,\tau_1 \in CP_1$
and $\lambda \in (0,1)$. By Proposition 2.7 we have $\tau_0(x)=V\pi(x)T'V^*$ for some non-negative $T' \in \pi(\cla)'$. $\tau_0$ being unital we also have $VT'V^*=I$ i.e. $V(I-T')V^*=0$. In case (b) is true then $T'=I$ and thus $\tau_0$  and so $\tau_1=\tau_0=\tau$. This proves (b) implies (a). For the converse $\Lambda \in \pi(\cla)'$ such that 
$V \Lambda V^*=0$. Same holds if we replace symmetric or anti-symmetric part of $\Lambda$. Thus it enough if prove (b) for self-adjoint $\Lambda$. Since it is a bounded operator, we assume without loss of generality that $-I \le \Lambda \le I$. So both $0 \le I-\Lambda \le 2I$ and $0 \le I+\Lambda \le 2I$ are operators in $\pi(\cla)'$. We write $I-\Lambda=W^*W$ and $I-\Lambda=W'^*W'$ for some $W,W'$ in $\pi(\cla)'$ and check that
$\tau={1 \over 2} (\tau_W + \tau_{W'})$ where
$\tau_W =VW\pi(x)W^*V^*$ and $\tau_{W'}=VW'\pi(x)W^{'*}V^*$ 
Since $\tau$ is extremal, we conclude that $\tau_W=\tau_{W'}=\tau$. Thus $(\clk,\pi,V^*)$ 
and $(\clk,\pi,W^*V^*)$ are two Stinespring minimum triplet. Thus $W^*V^*=U^*V^*$ and 
$U\pi(x)U^*=\pi(x)$ for all $x \in \cla$ where $U$ is an unitary operator on $\clk$. 
So we get $(UW^*-I)V^*=0$. Since both $U,W$ commutes with $\pi(\cla)$, we get
$(UW^*-I)\pi(x)V^*f=0$ for all $f \in \clh$ and $x \in \cla$. By cyclic property 
$\{ \pi(x)V^*f:x \in \cla,f \in \clh \}$ is total in $\clk$, thus $UW^*-I=0$ i.e. 
$W$ is unitary. So $\Lambda=I-W^*W=0$. 
\end{proof} 

\section{ Operator systems of unital normal complete positive maps: }

\vsp 
\begin{pro} 
Let $\tau:\clm \raro \clb(\clh)$ be a unital normal completely positive map on a von-Neumann algebra $\clm$ acting on a Hilbert space $\clh$. Then there exists a family of linearly independent operators $\{v_{\alpha}:\alpha \in \cli_{\tau}\}$ over the coefficients in $\clm'$, commutant of $\clm$ such that 
$\sum_{\alpha \in \cli_{\tau} }v_{\alpha}v_{\alpha}^*=I$ and 
$\tau(x)=\sum_{\alpha}v_{\alpha}xv_{\alpha}^*$ for all $x \in \clm$. If $\{v'_{\beta}:\beta \in \clj_{\tau} \}$ is another such a family of operators representing $\tau$ then there exists a family of elements $\{w^{\alpha}_{\beta} \in \clm':\alpha \in \cli_{\tau},\;\beta \in \clj_{\tau} \}$ 
such that $v'_{\beta}=\sum_{\alpha}v_{\alpha}w^{\alpha}_{\beta}$ 
and $(w^{\alpha}_\beta)$ is a unitary operator on $\clh \otimes \clk_0$. Furthermore, $\tau$ is an extremal element in $CP_{\sigma}$ if and only if there exists no non-trivial solution with elements $\lambda^{\alpha}_{\beta} \in \clm'$ satisfying 
$\sum_{\alpha,\beta \in \cli_{\tau} } v_{\alpha}\lambda^{\alpha}_{\beta}v_{\beta}=0$. 

\vsp 
Furthermore, if $\tau:\clm \raro \clm$ and admits an inner representation i.e. with $v_{\alpha} \in \clm$ then $\tau$ is extremal in the set of completely positive unital map on $\clm$ if and only if there is no non-trivial solution $(\lambda^{\alpha}_{\beta})$ in the center of $\clm$.   
\end{pro}

\vsp 
\begin{proof} 
$\tau$ being normal and unital, $\pi$ is a normal unital representation of $\clm$ into $\clb(\clk)$ and thus we may write $\clk \equiv \clh \otimes \clk_0$ for some Hilbert space $\clk_0$ and $\pi(x) \equiv x \otimes I_{\clk_0}$. Thus we conclude the first part by fixing an orthonormal basis $e_{\alpha}$ for $\clk_0$ and defining 
$<f,v^*_{\alpha}g>=<f \otimes e_{\alpha}, V^* g >$. To show linear independence of $\{v^*_{\alpha}: \alpha \in \cli_{\tau} \}$, let $\sum_{\alpha}c_{\alpha} v^*_{\alpha}=0$ 
for some $\sum |c_{\alpha}|^2 < \infty$. Then we have 
$$<f \otimes \sum_{\alpha} c_{\alpha} e_{\alpha},(x \otimes I_{\clk_0})V^*g>$$
$$=<x^*f \otimes \sum_{\alpha} c_{\alpha} e_{\alpha}, V^* g>$$
$$=\sum_{\alpha} c_{\alpha} <x^*f,v^*_{\alpha}g>$$
$$=<x^*f, \sum_{\alpha} c_{\alpha}v^*_{\alpha}g>$$  
$$=0$$ 
i.e. $f \otimes \sum_{\alpha} e_{\alpha}$ is orthogonal to the total set of 
vectors $\{(x \otimes I_{\clk_0}) V^*g: g \in \clh,x \in \clm \}$ in $\clk$  
and so $f \otimes \sum_{\alpha}c_{\alpha}e_{\alpha}=0$. Thus $\sum_{\alpha}c_{\alpha}e_{\alpha}=0$ and $(e_{\alpha})$ being linearly independent 
we get $c_{\alpha}=0$. 

\vsp 
A little modification of the proof works to show for the minimal representation, the family $v_{\alpha}$ is linearly independent over coefficients in $\clm'$ i.e. $\sum_{\alpha}c_{\alpha}v^*_{\alpha}=0$ with $c_{\alpha} \in \clm'$ and 
$\sum_{\alpha} c_{\alpha}^*c_{\alpha} \in \clm'$ if and only if $c_{\alpha}=0$. It is enough if we prove for $c_{\alpha} = 0$ for all $\alpha \in \cli_{\tau}$ except finitely many. To show linear independence of $\{v^*_{\alpha}: \alpha \in \cli_{\tau} \}$ over $\clm'$, let $\sum_{\alpha}c_{\alpha} v^*_{\alpha}=0$  
for some $c_{\alpha} \in \clm'$ with finitely many non zero elements. 
Then we have 
$$<\sum_{\alpha} c^*_{\alpha} f \otimes e_{\alpha},(x \otimes I_{\clk_0})V^*g>$$
$$=\sum_{\alpha} <x^*c^*_{\alpha}f \otimes e_{\alpha}, V^* g>$$
$$=\sum_{\alpha} <c^*_{\alpha}x^*f,v^*_{\alpha}g>$$
$$=<x^*f, \sum_{\alpha} c_{\alpha}v^*_{\alpha}g>$$  
$$=0$$ 
i.e. $\sum c^*_{\alpha}f \otimes e_{\alpha}$ is orthogonal to the total set of 
vectors $\{(x \otimes I_{\clk_0}) V^*g: g \in \clh,x \in \clm \}$ in $\clk$  
and so $\sum_{\alpha}c_{\alpha}^*f \otimes e_{\alpha}=0$. Thus 
$<\sum_{\alpha}c^*_{\alpha}f \otimes e_{\alpha}, g \otimes f_{\beta}>=0$ 
for all $\beta \in \cli_{\tau}$ and $g \in \clh$. The family $(e_{\alpha})$ being an orthonormal basis for $\clk_0$, we get $<c^*_{\alpha}f,g>=0$ for all $f,g \in \clh$. 
Thus we get $c_{\alpha}=0$. 

\vsp 
Choosing another such a family of elements $(v'_{\beta}: \beta \in \clj_{\tau})$ means that we are choosing another orthonormal basis $(e'^{\beta})$ for $\clk_0$ and thus there exists an unitary operator $W$ on $\clh \otimes \clk_0$ by (7) such that $W^*V^*=V'^*$. Since by our construction $W x \otimes I_{\clk_0} W^* = x \otimes I_{\clk_0}$ for all $x \in \clm$, we 
conclude that $W \in \clm' \otimes \clb(\clh_0)$. We set matrix elements $W^{\alpha}_{\beta} \in \clm'$ defined by 
$$<f,W^{\alpha}_{\beta}g>=<f \otimes e_{\alpha}, W e'_{\beta} \otimes g>$$
for all $f,g \in \clh$.   

\vsp 
We give now the proof of the last two statements. First statement is a simple corollary of Theorem 2.8 and normality of $\tau$. For the last statement, 
we compute with any unitary element $u \in \clm'$ by uniqueness of the Radon-Nykodym representation of $\eta$ for which  
$\eta \le c \tau$ on $\clm_+$ with a scaler $c > 0$, 
$\eta(x)=u\eta(x)u^*=u\sum_k v_k xt^k_jv_j^*u^*=\sum_k v_k x ut^k_ju^*v^*_j$
and thus by uniqueness of $t=(t^k_j)$ with entries in $\clm'$, we get $ut^k_ju^*=t^k_j$ and thus matrix entries in $t=(t^k_j)$ are elements in $\clm \bigcap \clm'$. 
This completes the proof.
\end{proof} 

\vsp 
Given a unital $C^*$-algebra $\cla$, a subspace $\cls$ of $\cla$ is called {\it operator system } if $\cls$ is closed under involution $*$ and identity of $\cla$ denoted by $I \in \cls$. Let $\cls_1,\cls_2$ be two operator systems in $C^*$-algebras $\cla_1$ and $\cla_2$ respectively. A linear one to one and onto map $\cli: \cls_1 \raro \cls_2$ is called 
order isomorphism if $\cli$ and $\cli^{-1}$ are both non-negative i.e. taking non-negative elements to non-negative elements. It is called a complete order isomorphism if $\cli \otimes I_n: M_n(\cls_1) \raro M_n(\cls_2)$ is also an order-isomorphism for each $n \ge 1$.  

\vsp 
For a unital completely positive normal map $\tau:\clm \raro \clb(\clh)$ with minimal representation $\tau(x)=\sum_k v_kxv_k^*$ for all $x \in \clm$, we set notations $\clm_{\tau}$ and $\cls_{\tau}$ for the operator space and system 
$$\clm_{\tau} = \{ \sum \lambda^i v^*_i: \lambda^i \in \clm' \}$$
$$\cls_{\tau} = \{ \sum v_i \lambda^i_j v_j^*: \lambda^i_j \in \clm' \}$$
The operator space $\clm_{\tau}$ over the algebra $\clm'$ is independent of the choice that we make for minimal representation of $\tau$ and thus the dimension $d_{\tau}$ of
$\clm_{\tau}$ is an invariance for $\tau$. This integer $d_{\tau}$ is called now onwards 
Arveson index of $\tau$. The operator system $\cls_{\tau}$ is a two sided module of dimension $d_{\tau}^2$ over the algebra $\clm'$, where left and right natural actions of $\clm'$ on $\cls_{\tau}$ are given by 
$$y \circ \sum v_ixv_j^* = \sum v_i yx v_j^*$$ 
$$\sum v_ixv_j^* \circ y = \sum v_ixyv_j^*$$
One simple question that arises here: does operator space $\cls_{\tau}$ depends on 
the choice that we make for $(v_k)$ to represent $\tau$? The following proposition says that $\cls_{\tau}$ is independent of the representation we choose and a little more. 

\vsp 
\begin{pro} 
Let $\tau:\clm \raro \clb(\clh)$ be a unital completely positive normal map. Then $\cls_{\tau}$ is determined uniquely modulo an unitary operator $u \in \cla'$ i.e. 
if $\tau(x)=\sum_k v_k xv_k^*$ and $\tau(x)=\sum_k v'_kx(v'_k)^*$ be two minimal representations of $\tau$ then there exists an unitary 
operator $\Lambda=((\lambda^i_j)) \in \!M_n(\clm')$ so that $\Lambda^*V^*=V'^*$ and operator system $\cls_{\tau}$ is uniquely determined by $\tau$ by its minimal representation representation i.e. there exists an unitary matrix 
$\lambda=(\lambda^i_j)$ with entries in $\clm'$ of order equal to numerical index so that 
\be 
v_k^*= \sum_j \lambda^k_j(v'_j)^*
\ee 
In particular, group of unitary elements $u \in \clm'$ acts on $(v_i^*)$ 
\be 
u^*v_k^*u = \sum_j \lambda^k_j(u)v_j^*,
\ee 
where $u \raro \Lambda(u)=((\lambda^j_k(u))) \in U_{d_{\tau}}(\clm')$ satisfies the 
following cocycle relation
\be 
\Lambda(uv)= (v^* \otimes I_{d_{\tau}}) \Lambda(u) (v \otimes I_{d_{\tau}}) \Lambda(v)
\ee
  
\vsp 
Conversely, let $\tau(x)=\sum_k v_kxv_k^*$ be a minimal representation of an extremal element $\tau$ in $CP_{\sigma}$ and $\eta=\sum_k w_kxw_k^*$ be another minimal representation of an element $\eta \in CP_{\sigma}$ so that $(w^*_k)$ is in the linear span of $(v_k^*)$ over $\clm'$ i.e. $\clm_{\tau}=\clm_{\eta}$ then $\eta=\tau$.
\end{pro}  

\vsp 
\begin{proof} 
By the uniqueness part of Stinespring intertwining relation find an unitary operator $\Lambda$ on $\clh \otimes \clk_0$ so that  
relation we have $\Lambda(x \otimes I)\Lambda^* = x \otimes I$ and $\Lambda V'^*=V^*$. Thus we get in the basis $v_k^*=\sum_j \lambda^k_jv'^*_j$ where $\Lambda=(\lambda^j_k)$ 
unitary elements in $\!M_{d_{\tau}}(\clm')$ determined by $<f \otimes e_i,\Lambda g \otimes e_j>=<f,\lambda^i_jg>$ where $d_{\tau}$ is the cardinality of an orthonormal 
basis for $\clk_0$. The relation (13) follows from (12) and linear independence of the family $(v_k^*)$ over $\clm'$. 

\vsp 
Now we consider the relation $\Lambda V'^*V' \Lambda^*=V^*V$ where $\Lambda \in U_{d_{\tau}}(\clm')$. Thus the change of the basis for $\clk$ will not change 
the operator systems $\cls_{\tau}$ and $\cls_{\tau'}$. 

\vsp 
For the last part we fix $\lambda^i_j \in \clm'$ so that $W^*=\Lambda V^*$ and claim that $\Lambda$ is an isometry if $\tau$ is an extremal element. To that end
we check $I=WW^*= V \Lambda^* \Lambda V^*$ i.e. $V(I-\Lambda^*\Lambda)V^*=0$. $\tau$ being extremal we get $\Lambda^*\Lambda-I=0$ by Corollary 2.9. Since $\lambda^i_j \in \clm'$, 
we get by a simple computation that for all $x \in \clm$, 
$$\eta(x)=\sum w_kxw_k^*= W(x \otimes I)W^*$$
$$=V\Lambda (x \otimes I )\Lambda^*V^* = V (x \otimes I) V^*= \sum_k v_kxv_k^*=\tau(x)$$ 
i.e. $\eta=\tau$. 
\end{proof} 

\vsp 
Now we formulate a more deeper problem. Two unital UCP maps 
$\tau_1:\clm_1 \raro \clb(\clh_1)$ and $\tau_2:\clm_2 \raro \clb(\clh_2)$  
are called {\it cocycle conjugate} if there exist automorphisms $\alpha:\clm_1 \raro \clm_2$
and $\beta:\clm_2 \raro \clm_1$ such that 
\be 
\tau_2 \circ \alpha = \beta \circ \tau_1
\ee
on $\clm_1$. Two cocycle conjugate elements are called {\it conjugate} if (14) is valid with $\alpha = \beta^{-1}$. For two conjugate elements $\tau_1$ and $\tau_2$, the operator systems $\cls_{\tau_1}$ and $\cls_{\tau_2}$ are complete order isomorphic. For a proof
we consider the lift of cocycle conjugate relation (14) to universal enveloping von-Neumann algebra $\clm_1^{**} = \clm_2^{**}$ as $\clm_1$ and $\clm_2$ are isomorphic. In other words, we can assume with out loss of generality that $\clm_1=\clm_2=\clm$ and $\clm$ is 
in its standard form acting on $\clh$ and thus for some unitary operators $u,v$ on $\clh$, we have $\alpha(x)=uxu^*$ and $\beta(x)=vxv^*$ for all $x \in \clm$. Thus operator spaces are cocycle conjugated by $u,v$ i.e. $u\clm_{\tau_1}v^*=\clm_{\tau_2}$. Since $\cls_{\tau}= \{x^*y:x,y \in \clm_{\tau} \}$, the operator systems $\cls_{\tau_1}$ 
and $\cls_{\tau_2}$ are complete order isomorphic via the map $\sum_{i,j} v_i\lambda^i_jv_j^* \raro \sum_{i,j} (uv_iv^*) v\lambda^i_jv^*(vv_j^*u^*)$, where
$\tau_1(x)=\sum_j v_jxv_j^*$.      

\vsp 
However, the converse is false even in the simplest situation when $\clm=\!M_2(\IC)$ 
and a detail results are given in a recent paper [Mo4]. However, one additional natural
condition on involved operator system gives a positive result. This inverse problem is
further addressed in a forth coming paper [Mo5].

\section{Tomita coupling and extremal marginal states:} 

\vsp 
We now aim to describe extremal marginal states in a more general mathematical frame work then originally proposed in [Par] and subsequently followed in papers [PSa,Ru,Oh]. Let $\phi$ be a faithful normal state on a von-Neumann algebra $\clm$ acting on a complex separable Hilbert space $\clh$ with inner product $<.,.>$ taken linear in the second variable and without loss of generality we assume $\clm$ be in its standard form $(\clm,\clp,\clj,\Delta,\Omega)$ [BR] associated with $\Omega$ where $\Omega$ is a cyclic and separating vector for $\clm$ and $\Delta,\clj$ are modular and conjugate operators of Tomita associated with polar decomposition $S=\clj \Delta^{1 \over 2}$, which is the closer of of the densely defined anti-linear closable operator 
$S_0x\Omega \raro x^*\Omega$ and $\clp=\overline{ \{ x \clj x \clj \Omega: x \in \clm \} }$ is the self-dual pointed positive cone in $\clh$. 
Here we recall that analytic elements $\clm_a$ of the modular group $\sigma_t(x)=\Delta^{it}x\Delta^{-it}$ on $\clm$ is dense in weak$^*$ topology 
in $\clm$ and we have the following modular relation for any two $x,y \in \clm_a$ given by 
\be 
\phi(x^*y)=\phi(\sigma_{i \over 2}(y) \sigma_{-{i \over 2}}(x^*))
\ee  

\vsp 
We consider the algebraic tensor product $\clm \ul{\otimes} \clm$ and complete it with $C^*$-cross norm of the von-Neumann algebra $\clm \otimes \clm$. The set of states $\psi$ on C$^*$ tensor product of von-Neumann algebras $\clm \circ \clm$ acting on $\clh \otimes \clh$ such that its restrictions on $\clm \circ I$ and $I \circ \clm$ are equal to $\phi$. In short, we will call such a element in words a coupling state with marginal $\phi$ and denote the convex set by $C_{\phi}$. Simplest example is the product state and simplest non-product state is given by a Tomita state 
$$\psi(x \circ y)=<\clj x \clj \Omega, y \Omega>$$ 
by extending linearly and then to its norm closures $\clm \circ \clm$. 
Important difference that we note now that $\psi$ may not have a normal 
extension to von-Neumann completion of
$\clm \circ \clm$ i.e. to $\clm \otimes \clm = (\clm \circ \clm)''$. As an example 
we take $\clm=L^{\infty}[0,1]$ and check that indicator function of the diagonal set in $[0,1] \times [0,1]$ can be expressed as limit of decreasing projections $E_n \in \clm \circ \clm$ with $\psi(E_n)=1$ and $\bigcap_{n \ge 1} E_n= \{[x,x]:0 \le x \le 1 \}$. 
In spite of this odd feature, we have a simple but beautiful observation. 

\vsp 
\begin{lem} 
A state $\psi \in C_{\phi}$ if and only if there exists a unital normal map $\tau_{\psi}:\clm \raro \clm$ preserving $\phi$ 
so that 
\be 
\psi(x \circ y)=<\clj x \clj \Omega, \tau_{\psi}(y) \Omega>
\ee 
for all $x,y \in \clm$. Further the map $\psi \raro \tau_{\psi}$ is an one to one and onto affine map between two convex set $C_{\phi}$ and 
$$CP_{\phi}=\{ \tau: \clm \raro \clm,
\mbox{ completely positive unital map and } \phi \circ \tau = \phi \}$$   
\end{lem}

\vsp 
\begin{proof} 
We fix $y \ge 0$ so that $\phi(y)=1$ and consider that state $\phi_y: x \raro \psi(x \otimes y)$ and note that $\phi_y(x) \le ||y|| \phi(x)$ for $x \ge 0$  
as $0 \le y \le ||y|| I$. Since by our assumption $\phi$ is normal, $\phi_y$ is also normal. Thus by Dixmier lemma we conclude that 
$$\phi_y(x)=< y'^* \Omega, x \Omega>$$ 
for some non-negative element $y'^* \in \clm'$. Thus we get 
$$\psi(x \circ y)= <\clj x \clj \Omega, \clj y'^* \clj \Omega>$$ 
We set $\tau(y)=\clj y'^* \clj \in \clm$ to arrive at (16). Since $\tau(y)$ is determined uniquely by the separating property of $\Omega$ for $\clm$, we may extend the map for an arbitrary element $y$ using linearity by writing it as a linear combination of four non
negative elements in $\clm$. That $y \raro \tau(y)$ is a normal map follows by invariance of the faithful normal state $\phi$ for $\tau$ and positivity of the map as follows: for an increasing net $y_{\alpha}$ with least upper bound $y$, $\tau(y_{\alpha})$ is also an increasing net with $\tau(y)$ as an upper bound and thus has a least upper bound say $z$. Then $\phi(z-\tau(y))= \mbox{l.u.b.}_{\alpha}\phi(\tau(y_{\alpha})-\tau(y))$ by normality of $\phi$ and thus by invariance equal to $\mbox{l.u.b.}_{\alpha}\phi(y_{\alpha}-y)$ which is $0$ by normality of the state. For complete positive property of the map $\tau$, we first check that 
$$<\Omega, \sum_{i,j}(z'_i)^*\tau(y^*_jy_i)z'_j \Omega>$$
$$=\sum_{i,j} <\Omega, (z'_i)^*z'_j\tau(y_jy^*_i) \Omega>$$
$$=\sum_{i,j} <(z'_j)^*z'_i\Omega, \tau(y_jy^*_i) \Omega>$$
$$=\sum_{i,j} < \clj z_j^*z_i \clj \Omega, \tau(y_jy^*_i)\Omega>$$
$$=\sum_{i,j} \psi(z_j^*z_i \circ y_jy_i^*)$$ 
$$=\psi (X^*X) \ge 0 $$ 
where $X = \sum_i z_i \circ y^*_i $  for all $1 \le i \le n$, $y_i \in \clm, z'_i \in \clm'$ and $z_i=\clj z'_i \clj$. Since $\Omega$ is cyclic for $\clm'$, we conclude that the operator $((\tau(y_jy_i^*)))$ is a non-negative element in $\!M_n(\clm)$. Thus $\tau$ is $n-$positive for each $n \ge 1$. Rest of the statements are now obvious. 
\end{proof} 

\vsp 
\begin{thm} 
The affine map $\psi \raro \tau_{\psi}$ defined in Lemma 3.1 takes extremal elements of the convex set $C_{\phi}$ to extremal elements of $CP_{\phi}$. Further

\NI (a) $C_{\phi}$ is a closed convex subset in the weak topology of $\clm \circ \clm$;

\NI (b) $CP_{\phi}$ is also compact once equipped with the point-wise topology inherited from  
$\sigma$-weak operator topology of $\clm$ (i.e. we say a net $\tau_{\alpha} \in CP$ 
converges to $\tau \in CP$ if the net $\tau_{\alpha}(x)$ converges to $\tau(x)$ in 
$\sigma-$weak operator topology for each $x \in \clm$, where $CP$ denotes the set of 
unital completely positive maps on $\clm$ ).    
\end{thm} 

\vsp 
\begin{proof} 
First part is obvious as the map is one to one and onto. Crucial observation that we make here for a net $\tau_{\alpha}$ 
in $CP_{\phi}$ such that $\tau_{\alpha}(x) \raro \tau(x)$ in $\sigma-$weak operator topology for some $\tau \in CP$. Then $\tau$ is also 
unital and $\phi$ preserving and $\tau$ is normal as $\phi$ is faithful. Thus $\tau \in C_{\phi}$. 

\vsp 
We consider a net of states $\psi_{\alpha} \in C_{\phi}$ so that $\psi_{\alpha} \raro \psi$ 
in weak topology on $\clm \circ \clm$ and since each $\psi_{\alpha}$ preserves marginals 
so is their limit. Thus $\psi \in C_{\phi}$. Now we also check that  
\be 
\psi_{\alpha}(x \circ y)=\phi(\clj x \clj \Omega, \tau_{\alpha}(y)\Omega> \raro \psi(x \circ y)
\ee 
for all $x,y \in \clm$, where the state $\psi$ on $\clm \circ \clm$ defined by 
$$\psi(x \circ y)=<\clj x \clj \Omega,\tau(y)\Omega>$$ 
determines a unique element $\tau \in CP_{\phi}$ by Lemma 3.1. These shows that $\tau_{\alpha}(x) \raro \tau(x)$ in weak operator topology. Since the family $(\tau_{\alpha})$ is uniformly bounded and $\Omega$ is cyclic and separating for $\clm$, the limit (17) says that 
$$\sum_k <f_k,\tau_{\alpha}(x)g_k> \raro \sum_k <f_k,\tau(x)g_k>$$ 
as the net $\alpha \raro \infty$ for $\sum_k ||f_k||^2 < \infty$ and $\sum_k ||g_k||^2 < \infty$ by dominated convergence theorem. Thus $\tau_{\alpha}(x) \raro \tau(x)$ 
in $\sigma-$weak operator topology. In other words, convergence holds in Bounded Weak 
topology of Arveson.  

\vsp 
Conversely for a given net $\tau_{\alpha}$ in $CP_{\phi}$ which converges to $\tau$ in point-wise $\sigma-$weak operator topology, $\tau \in CP_{\phi}$ 
by our remark at the beginning of the proof. We check that the associated elements $\psi_{\alpha}$ also converges to $\psi$ in the weak topology of 
$\clm \circ \clm$ first on a norm dense subspace of algebraic tensor product of $\clm \underline{\otimes} \clm$ and then for any arbitrary elements by a standard argument as the net is uniformly bounded.    

\vsp 
The convex set of states on $\clm \circ \clm$ is compact in weak$^*$ topology and thus $C_{\phi}$ being a closed subset of the compact set in the weak$^*$ topology, $C_{\phi}$ is also compact. That $CP_{\phi}$ is compact also follows from compactness of $C_{\phi}$ via the above correspondence which respect the topologies.  
\end{proof} 

\vsp 
One can as well use BW (Bounded Weak) topology [Ar,Pa, Chapter 7] of Arveson directly to give a direct proof that $CP_{\phi}$ is a closed subset of the unit ball of bounded linear maps on $\clm$. Any limit points of a convergent net of $\phi$-invariant unital completely positive normal maps $\tau_{\alpha}$ will be also completely positive and $\phi$-invariant. $\phi$ being faithful and normal, any positive $\phi$ invariant map is also normal. However such an argument will not be valid with BW topology for a more general situation, where $\phi$ is just a normal $\sigma-$finite weight [BR,Ta]. It is not clear what would be a possible truncation method along the classic work [Ke,Ko]. 

\vsp 
By our last theorem we conclude that extremal elements in $CP_{\phi}$ exists and any element in $CP_{\phi}$ admits Krein-Milman property. In the following text we aim to find a criteria for an extremal element $\tau \in CP_{\phi}$. To that end, we recall KMS-adjoint completely positive map [OP].  

\vsp 
\begin{pro} 
Let $\phi$ be a faithful normal state on a von-Neumann algebra $\clm$. Then $\tau$ is a positive normal map on $\clm$ such that $\phi \tau \le \phi$ if and only if exists a unique normal positive map $\tilde{\tau}$ on $\clm$ satisfying 
the following duality relation
\be 
\phi(\tau(x)\sigma_{-{i \over 2}}(y))=\phi(\sigma_{i \over 2}(x)\tilde{\tau}(y))
\ee 
where $x,y \in \clm_a$, the $*$-algebra of analytic elements for modular group $\sigma=(\sigma_t:t \in \IR)$ associated with $\phi$ and $\tilde{\tau}(I) \le I$. Moreover, $\phi \tilde{\tau} \le \phi$ if and only if $\tau(I) \le I$. 

\vsp 
Furthermore, $\tau$ is completely positive if, and only if $\tilde{\tau}$ is completely positive. In such a case, the numerical indices of $\tau$ and $\tilde{\tau}$ are equal. 
\end{pro}  

\vsp 
\begin{proof} 
For the first part of the statement we refer to chapter 8 of monograph [OP]. For the second part we re-investigate the proof of Stinespring representation 
with our special situation. With out loss of generality we assume that $\clm$ is in the standard form associated with $\phi$. i.e. $\phi(x)=(\Omega,x \Omega)$, where
$\Omega$ is a cyclic and separating vector for $\clm$. We set a kernel on 
$\clm_a \otimes \clm_a$ defined by
$$k(x \otimes z, y \otimes w)=<\Omega,x^*\tau(z^*w)y\Omega>$$
That Hilbert space completion of the kernel is same as that of Stinespring follows by cyclic property of $\Omega$. Similarly we set
kernel $\tilde{k}$ associated with $\tilde{\tau}$ by 
$$\tilde{k}(x \otimes z, y \otimes w)=<\Omega,z^*\tau(x^*y)w\Omega>$$
in reverse direction and check by KMS relation that 
$$<\Omega,x^*\tau(z^*w)y\Omega>=<\Omega,\tilde{w}^*\tilde{\tau}(\tilde{y}^*\tilde{x})\tilde{z}\Omega>$$   
where for any analytic element $x \in \clm$ we write $\tilde{x}= \sigma_{-{i \over 2}}(x^*)$. Such a relation is used already in [Mo1,Mo2,Mo3].  
This clearly shows that Stinespring Hilbert spaces $\clk$ and $\tilde{\clk}$ associated with kernels $k$ and $\tilde{k}$ respectively are conjugated by an anti-unitary operator defined by 
$U: x \otimes y \raro \tilde{y} \otimes \tilde{x}$. Since anti-unitary operator $U$ inter-twins the Stinespring representations $(\clk,\pi,V)$ and $(\tilde{\clk},\tilde{\pi},\tilde{V})$ it also inter-twins $\clk_0$ and $\tilde{\clk}_0$ and thus we get dimension of $\clk_0$ and $\tilde{\clk}_0$ are same. Thus index of $\tau$ and $\tilde{\tau}$ are equal. 
\end{proof}

\vsp  
In the following, we give a criteria for an element $\tau$ in $CP_{\phi}$ to be extremal.
If $\clm$ is a matrix algebra and $\phi$ is the normalize trace, the criteria coincides with Landau-Streater criteria [LS] known in the literature for quite some time. Our proof follows quite a different method inspired by Proposition 4.1.               

\vsp 
\begin{thm} 
Let $\tau$ and $\tilde{\tau}$ be the elements in $CP_{\phi}$ of equal numerical 
index admitting the following minimal representations 
$$\tau(x)=\sum_{\alpha} v_{\alpha}xv_{\alpha}^*$$ 
and 
$$\tilde{\tau}(x) = \sum_{\alpha} \tilde{v}_{\alpha} x \tilde{v}_{\alpha}^*$$
for all $x \in \clm$, where $\tau$ and $\tilde{\tau}$ are related by the duality relation 
(16) in Proposition 5.3. 

\vsp 
Then $\tau$ is an extremal element in $CP_{\phi}$ if and only if
there exists no non-trivial $\lambda=(\lambda^k_j)$ with entries in $\clm'$ satisfying 
the relation 
\be 
\sum_{\alpha,\beta} v_{\alpha} \lambda^{\alpha}_{\beta} v^*_{\beta}=0,\;\;\sum_{\alpha,\beta} \tilde{v}_{\alpha} \tilde{\lambda}^{\alpha}_{\beta} \tilde{v}^*_{\beta}=0,
\ee
where $\lambda \raro \tilde{\lambda}$ is an unital order-isomorphism map on $\!M_{n_{\tau}}(\clm'))$.

\end{thm} 

\vsp 
\begin{proof} 
Let $\eta$ be an element in $CP$ such that $\eta \le k \tau$ for some $k \ge 0$. Then $\eta_0=k \tau- \eta$ is also a positive map and $\eta_0 \le k \tau$. The duality being 
an affine map, we get $\tilde{\eta} + \tilde{\eta}_0 = k \tilde{\tau}$, where $\tilde{\eta}_0$ is a positive map by Proposition 5.3. Thus we also have 
$\tilde{\eta} \le k \tilde{\tau}$.   

\vsp 
Then $\eta$ is a normal map and there exists a unique non-negative element 
$t=(t^{\alpha}_{\beta})$ with entries in $\clm'$ such that 
\be 
\eta(x)=\sum_{\alpha,\beta} v_{\alpha} xt^{\alpha}_{\beta}v_{\beta}^*
\ee
for all $x \in \clm$. Conversely. a non-negative bounded operator  $T=(t^{\alpha}_{\beta})$ with entries in $\clm$ determines a normal map 
$\eta$ by (20) for some $k \ge 0$. Thus there is a one to one relation between completely positive map satisfying $\eta \le k \tau$ and bounded non-negative elements $T=(t^{\alpha}_\beta)$ with entries in $\clm'$ determined by (20).      

\vsp 
Similarly, for a give $\eta \le k \tau$, there exists a unique $\tilde{T}=(\tilde{t}^{\alpha}_{\beta})$ with entries in $\clm'$ such that
\be 
\tilde{\eta}(x) = \sum_{\alpha,\beta} \tilde{v}_{\alpha}x\tilde{t}^{\alpha}_{\beta}v_{\beta}^*
\ee 
for all $x \in \clm$. 

\vsp 
This shows that $T \raro \tilde{T}$ is a well defined affine map on the non negative elements of $\!M_{d_{\tau}}(\clm')$ determined by (20) and (21) for fixed minimal representations of $\tau$ and $\tilde{\tau}$. 

\vsp 
We extends the map $T \raro \tilde{T}$ to self-adjoint elements $T=T^1_+ - T^2_+$ with $T^1_+,T^2_+ \ge 0$ to $\tilde{T} = \tilde{T^1_+} - \tilde{T^2_+}$. To show that the map is well defined, 
let $T=S^1_+-S^2_+$ be another expression with $S^1_+,S^2_+ \ge 0$. Then $T^1_+ +S^2_+ = S^1_+ + T^2_+$. Since $T \raro \tilde{T}$ is an affine map on $\!M_{n_{\tau}}(\clm')_+$,
we get $\tilde{T^1_+} + \tilde{S^2_+} = \tilde{S^1_+} + \tilde{T^2_+}$. Rearranging the terms, we conclude the map $T \raro \tilde{T}$ is well defined on self-adjoint elements of $\!M_{n_{\tau}}(\clm')$. That it is an injective map, follows readily as $\tilde{T}=0$ implies $\tilde{T_+} = \tilde{T_-}$. The map being injective on non-negative elements, we have $T_+=T_-$. Thus $T+^2=T_-^2=T_+T_-=0$ i.e. $T=0$. The map is clearly onto on self adjoint elements. We extend now by linearity to all element from 
$\!M_{n_{\tau}}(\clm')$ to itself.  

\vsp 
Since $\tilde{T} \raro T$ is also an injective map on non-negative matrices, there exists 
an injective and onto extension of this map from self adjoint elements to self adjoint elements. Thus we conclude that $T \raro \tilde{T}$ extends to an order isomorphic map uniquely on $\!M_{n_{\tau}}(\clm')$. 

\vsp 
Thus $\lambda^{\alpha}_{\beta}=t^{\alpha}_{\beta}-\delta^{\alpha}_{\beta}I$ is a solution to (19) and it is non-trivial if and only if $\tau$ is not an extremal element in $CP_{\phi}$. This completes the proof. 
\end{proof} 

\vsp 
By a theorem of Kadison [Ka1], the map $T \raro \tilde{T}$ is a 
disjoint sum of a morphism and anti-morphism determined by a projection $\tilde{e}$ in the centre of $\clm$ i.e. there exists a projection $\tilde{e}$ in the centre of $\clm$ such that 
$$T  \raro \tilde{T} \tilde{e} \otimes I_{n_{\tau}}$$ 
is a morphism and 
$T \raro \tilde{T} (I-\tilde{e}) \otimes I_{n_{\tau}}$ 
is an anti-morphism. In case $\clm$ is a factor, then $T \raro \tilde{T}$ is either 
an isomorphism or an anti-isomorphism. 

\vsp 
If $\tau$ admits an inner minimal representation 
\be 
\tau(x)=\sum_{\alpha \in \clc} v_{\alpha}xv_{\alpha}^*,\;x \in \clm
\ee 
i.e. with elements $v_{\alpha} \in \clm$ then $\tilde{\tau}$ also admits an inner minimal representation given by  
\be 
\tilde{\tau}(y)=\sum_{\alpha \in \clc}\tilde{v}_{\alpha}y\tilde{v}_{\alpha}^*,\;y \in \clm
\ee 
where $\tilde{v}_{\alpha} \in \clm$ is defined as the bounded operator extending the densely defined operator 
$f \raro \Delta^{1 \over 2} v_{\alpha}^* \Delta^{-{1 \over 2}}f$ for all $f \in \clm \Omega$. For the last part, we refer section 7 and appendix given in [BJKW].

\vsp 
In particular, if $\eta \le k \tau$ for some $k \ge 0$ then the representations of $\eta$ and $\tilde{\eta}$ given in (20) and (22) are inner with matrices 
$(t^{\alpha}_{\beta})$ and $(\tilde{t}^{\alpha}_{\beta})$ are with entries 
in the centre of $\clm$. Thus the map $T \raro \tilde{T}$ is an order isomorphism on 
$\!M_{d_{\tau}}(\clz)$, where $\clz = \clm \bigcap \clm'$.  

\vsp 
A simple computation using modular relation (15) now also leads to the 
minimal representation 
$$\tilde{\eta}(x)=\sum \tilde{v}_{\alpha} x t^{\alpha}_{\beta}\tilde{v}^*_{\beta}$$
since modular group acts trivially on the center $\clz$ of $\clm$. By the uniqueness part of representation (21), we conclude that 
$$\tilde{t}^{\alpha}_{\beta}=t^{\beta}_{\alpha}$$

\vsp 
We sum up now our results in the following corollary.

\vsp 
\begin{cor} 
Let $\tau$ be an element in $CP_{\phi}$ with inner representation given by (22). Then its dual element $\tilde{\tau}$ also admits a inner representation given by (23). Furthermore, $\tau$ is an extremal element in $CP_{\phi}$ if, and only if there exists no non-trivial $\lambda=(\lambda^k_j)$ with entries in the centre of $\clm$ satisfying 
the relation 
\be 
\sum_{\alpha,\beta} v_{\alpha} \lambda^{\alpha}_{\beta} v^*_{\beta}=0,\;\;\sum_{\alpha,\beta} \tilde{v}_{\alpha} \lambda_{\alpha}^{\beta} \tilde{v}^*_{\beta}=0
\ee
\end{cor} 

\vsp 
Krein-Milman theorem says that $\tau= \int \eta d\mu(\eta)$ for some probability measure $\mu$ on the closer of extreme points $CP^e_{\phi}$ of $CP_{\phi}$ in BW topology [Pa, Chapter 7], where $\mu$ may not be uniquely determined. A valid question that rises here when can we guarantee $\mu$ to have support on $CP^e_{\phi}$ only? The topology on the states of $\clm \circ \clm$ is not metrizable unless $\clm$ is separable as a $C^*$ algebra. Such a restriction makes our choice for $\clm$ rather limited since $\clm$ 
is a von-Neumann algebra. When $\clm=l^{\infty}(\IN)$, some clever truncated methods are used [Ke, Ko] to show the support of $\mu$ is confined to $CP^e_{\phi_0}(l^{\infty}(\IN))$, where $\phi_0$ is the counting measure on $\IN$.

\section{Pure marginal states:} 

\vsp 
In the following, we answer a question raised in recent papers [Par,PSa,Oh]. 

\vsp 
\begin{thm} 
Let $\clm=\!M_n(\IC)$ and $\phi$ be a faithful normal state on $\clm$. 
An extremal element $\psi$ in $C_{\phi}$ is a pure state if and only if $\tau_{\psi}$ is 
an automorphism on $\clm$.   
\end{thm} 

\vsp 
\begin{proof} 
We may follow proof of Theorem 1.1 in [PSa] with obvious modification to prove: $\psi \in C_{\phi}$ is a pure state of $\clm \otimes \clm$ if and only if there are orthonormal bases $(f_i:1 \le i \le n)$ and $(g_i:1 \le i \le n)$ of $\IC^n$ such that 
$\psi(X)=<\zeta_{\psi},X \zeta_{\psi}>$ for all $X \in \clm \otimes \clm$, where $\zeta_{\psi}= \sum_i \lambda_i^{1 \over 2} f_i \otimes g_i$, $\phi(x)=tr(\rho x),\;x \in \clm$ and $\rho = \sum_{1 \le i \le n} \lambda_i |f_i><f_i| = \sum_{1 \le i \le n}\lambda_i |g_i ><g_i|$.     

\vsp 
For two such pure states $\psi$ and $\psi'$, we find unitary operators $u,v$ on $\IC^n$ which takes bases $(f_i) \raro (f'_i)$ and $(g_i) \raro (g'_i)$ respectively. We claim that $$\tau_{\psi'}(x)= u\tau_{\psi}(vxv^*)u^*$$
For a proof we use (14) and $u \otimes v \zeta_{\psi}=\zeta_{\psi'}$ to compute 
$$<\clj x \clj \Omega,\tau_{\psi'}(y)\Omega>$$
$$=\psi'(x \otimes y)$$ 
$$=\psi(u \otimes v (x \otimes y) u^* \otimes v^*)$$
$$=\psi(uxu^* \otimes vyv^*)$$
$$=<\clj ux u^* \clj \Omega,\tau_{\psi}(vyv^*)\Omega>$$
$$=<\clj x \clj u^*\Omega, u^*\tau_{\psi}(vyv^*)\Omega>$$
(since $\phi(uxu^*)=\phi(x),x \in \clm$, $u$ commutes with Tomita's modular operator and conjugate operator and $u\Omega=\Omega$ ).     
$$=<\clj x \clj \Omega, u^*\tau_{\psi}(vyv^*)u\Omega>$$
Thus we have $\tau_{\psi'}(y)=u^*\tau_{\psi}(vyv^*)u$ for all $y \in \clm$. Since 
$\psi(x \otimes y) =<\clj x \clj \Omega,y \Omega>$ is a pure state on 
$\clm \otimes \clm$, we get the required result. 
\end{proof} 

\vsp 
\begin{thm} 
Let $\phi$ be the normalized trace on $\clm=\!M_n(\IC)$ and $\psi$ be a state on $\clm \otimes \clm$ with marginal states $\phi$ i.e. $\psi \in C_{\phi}$. If $\psi$ is an extremal element in $C_{\phi}$ then $\psi$ is a factor state of $\clm 
\otimes \clm$. 
\end{thm} 

\vsp 
\begin{proof} 

\vsp 
A state $\psi$, $\psi(X)=tr(Xh_{\psi}),\;X \in \clm \otimes \clm$ for some non-negative density operator $h_{\psi} \in \clm \otimes \clm$, is an element in $C_{\phi}$ if and only if $\IE_1(h_{\psi})={1 \over d}I_d$ and $\IE_2(h_{\psi})={1 \over d} I_d$, where $\IE_1,\IE_2$ are conditional expectation with respect to the normalized trace from $\clm \otimes \clm$ onto $\clm \otimes I_n$ and $I_n \otimes \clm$ respectively. 

\vsp 
Let $\psi$ be an extremal element in $C_{\phi}$ and $(\clh_{\psi},\pi_{\psi},\zeta_{\psi})$ be its GNS representation. We claim that $\psi$ is a factor state i.e. 
the centre of $\pi_{\psi}(\clm \otimes \clm)''$ is trivial. Suppose it is not. Let $E$ be 
a non trivial projection in the centre of $\pi_{\psi}(\clm \otimes \clm)''$ and  
$\psi = \lambda \phi_1 + (1-\lambda)\phi_0$, where $\psi_1,\psi_0$ are states on  
$\clm \otimes \clm$ defined by 
$$\lambda \psi_1(X)=\psi_E(X)=<\zeta_{\psi},XE\zeta_{\psi}>=tr(h_{\psi}XE)$$
and
$$(1-\lambda) \psi_0(X)=\psi_{I-E}(X)=<\zeta_{\psi},X(I-E)\zeta_{\psi}>=tr(h_{\psi}X(I-E))$$
for $\lambda = \psi(E) \in (0,1)$ and $h_0,h_1 \in \clm \otimes \clm$ since $E$ is
an element in $\clm \otimes \clm$ as well.
Thus we have 
$$h= \lambda h_1 + (1-\lambda)h_0,$$
where $\psi_k(X) = tr(h_kX)$ for $k=0,1$. 
So ${1 \over d} I_d=\IE_1(h_{\psi})= \lambda \IE_1(h_1) + (1-\lambda)\IE_1(h_0)$, where we have used $h_0,h_1 \in \clm \otimes \clm$. The operator ${1 \over d}I_d$ being an extremal element in the set of non-negative definite matrices of trace $1$, we get $\IE_1(h_1)= \IE_1(h_0)={1 \over n}I_d$. Similarly, we also have $\IE_2(h_1)=\IE_2(h_0)={1 \over d}I_d$. This shows in particular, $\psi_0,\psi_1 \in C_{\phi}$ and $\psi = \lambda \psi_1 +(1-\lambda)\psi_0$ for some $\lambda \in (0,1).$ This brings a contradiction to extremal 
property of $\psi$ in our hypothesis at the beginning.  

\end{proof}

\vsp 
The converse statement of Theorem 6.2 is obviously false. For a counter example, 
we consider the product state $\psi(x \otimes y)= \phi(x) \phi(y)$ on $\clm \otimes \clm$, where $\phi$ is a faithful state on $\clm$. In such a case, $\tau_{\psi}=\phi$ 
and $\phi$ is a factor state but need not be an extremal element unless $\clm$ is 
itself $\IC$ i.e. $n=1$ in Theorem 6.2.

\vsp 
However, for a state $\psi$ in $C_{\phi}$ and we may consider the extremal decomposition of $\tau_{\psi}$ in $CP_{\phi}$ given by  
$$\tau_{\psi} = \int_{CP_{\phi}^e} \tau_e d\mu_{\tau}(e)$$
for some regular Borel probability measure $\mu_{\tau}$ on the set of extremal elements. By what we have proved above, each $\tau_e$ gives a unique factor state $\psi_e$ of $\clm \otimes \clm$ 
and so that 
\be
\psi = \int_{CP^e_{\phi}} \psi_e d\mu_{\tau}(e)
\ee
If $\psi$ is a factor state but not an extremal element, then 
support of $\mu_{\tau}$ is not atomic and the decomposition (25) is not central i.e. associated GNS representations is not a central decomposition [BRI].   

\vsp 
Theorem 6.1 and Theorem 6.2 give raises the following interesting question. 
Let $\psi$ and $\psi'$ be two states of $\clm \otimes \clm$ in $C_{\phi}$ with density matrices $h_{\psi}$ and $h_{\psi'}$ respectively. Then $\tau_{\psi}$ and $\tau_{\psi'}$  are cocycle conjugate if, and only if $h_{\psi}$ and $h_{\psi'}$ are unitary equivalent by an unitary matrix of the form $u \otimes v$. So a classification of extremal elements is equivalent to classify all density matrices $h$ in $\clm \otimes \clm$ for which the state
$\psi_h(x)=tr(xh)$ is atleast a factor state and $\IE_1(h)={1 \over d}I_d$ and $\IE_2(h)=
{1 \over d}I_d$. One necessary condition is equality of ranks of $h_{\psi}$ and $h_{\psi'}$. Theorem 6.1 suggests rank of $h_{\psi}$ is possibly a complete invariance for an extremal element $\tau_{\psi}$ in $CP_{\phi}$. We defer a possible answer of this finer question to [Mo5] which takes some hint from recent results proved in [Mo4]. 

\vsp 
Any faithful state $\phi$ on $\clm \subseteq \clb(\clh)$ admits a representation 
$$\phi(x) = \sum_k \lambda_k <f_k,xf_k>$$
for all $x \in \clm$ with some $\lambda_k > 0$ for all $k \ge 1$,
where $\clh$ is assumed to be a separable Hilbert space. If $\clm=\clb(\clh)$, an obvious modification of the argument used in the proof of Theorem 6.1, also proves that any extremal element $\psi$ in $C_{\phi}$ is pure if, and only if $\tau_{\psi}$ is an automorphism on $\clm$. Same questions for an arbitrary von-Neumann algebra $\clm$ is rather delicate since a state $\psi \in C_{\phi}$ may not have a normal extension to a state of $\clm \otimes \clm$.

\bigskip
{\centerline {\bf REFERENCES}}

\begin{itemize} 

\bigskip
\item{[Ar]} Arveson, W.: Sub-algebras of $C^*$-algebras, Acta Math. 123, 141-224, 1969. 

\item {[BR]} Bratteli, Ola., Robinson, D.W. : Operator algebras and quantum statistical mechanics I,II, Springer 1981.

\item{[BJKW]} Bratteli, Ola,; Jorgensen, Palle E.T.; Kishimoto, Akitaka and
Werner Reinhard F.: Pure states on $\clo_d$, J.Operator Theory 43 (2000),
no-1, 97-143.    

\item {[Ch]} Choi, M.D.: Completely positive linear maps on complex matrices, Linear Algebra and App (10) 1975 285-290. 

\item{[HMP]} Hopenwasser, Alan; Moore, Robert L.; Paulsen, V. I. : $C^*$-extreme points, Trans. Amer. Math. Soc. 266 (1981), no. 1, 291-307. 

\item {[Ka1]} Kadison, Richard V.: Isometries of operator algebras, Ann. Math. 54(2)(1951) 325-338.  

\item {[Ka2]} Kadison, Richard V.: A generalized Schwarz inequality and algebraic invariants for operator algebras,  Ann. of Math. (2)  56, 494-503 (1952). 

\item {[Ke]} Kendall, D.G. : On infinite doubly stochastic matrices and Birkhoff's problem, III. London Math. Soc. J., 35 (1960):81-84. 

\item {[Ko]} K\"{o}nig, D. : The theory of finite and infinite graphs, T\"{a}ubner 1936, Birkh\"{a}ser, Boston, 1990, p-327. 

\item {[LS]} Landau, L.J., Streater, R.F.: On Birkhoff theorem for doubly stochastic completely positive maps of matrix algebras, 
Linear Algebra and its Applications, Vol 193, 1993, 107-127 

\item {[MW]} Mendl, Christian B., Wolf, Michael M.: Unital quantum channel's convex structure and revivals of Birkhoff's theorem. Comm. Math. Phys. 289 (2009), no. 3, 1057--1086. 

\item {[Mo1]} Mohari, A.: Pure inductive limit state and Kolmogorov property. II 
Journal of Operator Theory. vol 72, issue 2, 387-404.   
    
\item {[Mo2]} Mohari, A.: Translation invariant pure state on $\otimes_{\IZ}\!M_d(\IC)$ and Haag duality, Complex Anal. Oper. Theory 8 (2014), no. 3, 745-789.

\item{[Mo3]} Mohari, A.: Translation invariant pure state on $\clb=\otimes_{\IZ}M_d(\IC)$ and its split property, J. Math. Phys. 56, 061701 (2015).

\item{[Mo4]} Mohari, A.: Hann-Banach-Arveson extension theorem and Kadison isomorphism, arXiv:1304.6849 (2015). 

\item{[Mo5]} Mohari, A.: G. Birkhoff problem in irreversible quantum dynamics, in preparation (2015). 

\item{[Oh]} Ohno, H.: Maximal rank of extremal marginal tracial states.  J. Math. Phys.  51  (2010),  no. 9, 092101, 9 pp. 

\item {[OP]} Ohya, M., Petz, D.: Quantum entropy and its use, Text and monograph in physics, Springer-Verlag 1995.

\item {[Pa]} Paulsen, V.: Completely bounded maps and operator algebras, Cambridge Studies in Advance Mathematics 78, 
Cambridge University Press. 2002 

\item{[Pas]} Paschke, W. L. : Inner Product Modules Over $B^*$-Algebras, 
Transactions of the American Mathematical Society, Vol. 182 (1973), pp 443-468

\item{[Ph]} Phelps, Robert R.: Lectures on Choquet's Theorem, Lecture notes in Mathematics 1757, Springer 2001.

\item{[Par]} Parthasarathy, K.R.: Extremal quantum states in coupled states in coupled systems, Ann.Inst. H. Poincar\'{e}, 41, 257-268 (2005).  
 
\item{[PSa]} Price, G.L., Sakai, S.: Extremal marginal tracial states in couple systems, Operators and matrices, 1, 153-163 (2007). 

\item{[Ra]} Raginsky, M. Radon-Nikodym derivatives of quantum operations. J. Math. Phys. 44 (2003), no. 11, 5003–5020. 

\item{[Ru]} Rudolph, O.: On extremal quantum states of composite systems with fixed marginals, J. Math. Phys. 45, 4035-4041 (2004).

\item {[St]} Stinespring, W. F.: Positive functions on $C^*$ algebras, Proc. Amer. Math. Soc. 6 (1955) 211-216. 
 
\item {[St\o]} St\o rmer, E.: Positive linear maps of operator algebras. Acta Math. 110, 233-278 (1963).
 
\item{[Ta]} Takesaki, M. : Theory of Operator algebras I, Springer, 2001.

\item{[To1]} Tomiyama, J.: On the projection of norm one in W$^*$-algebras. Proc. Japan Acad. 33 1957 608-612.

\item{[To2]} Tomiyama, J.: On the projection of norm one in W$^*$-algebras. II. Tôhoku . Math J. (2) 10 1958 204-209. 

\item{[To3]} Tomiyama, J.: On the projection of norm one in W$^*$-algebras. III. Tôhoku Math. J. (2) 11 1959 125-129. (

\item{[Ts]} Tsui, S. K.: Completely positive module maps and completely positive extreme maps. Proc. Amer. Math. Soc. 124 (1996), no. 2, 437-445.

\end{itemize}

\end{document}